\documentclass[10pt,leqno]{article}

\usepackage{amsmath,amsfonts,amscd,amssymb,theorem}

\long\def\comment#1\endcomment{}

\makeatletter
\begingroup
\gdef\th@dotted{\normalfont\itshape
  \def\@begintheorem##1##2{%
        \item[\hskip\labelsep \theorem@headerfont ##1\ ##2.]}%
\def\@opargbegintheorem##1##2##3{%
   \item[\hskip\labelsep \theorem@headerfont ##1\ ##2\ (##3).]}}
\endgroup
\makeatother

\theoremstyle{dotted}

\newtheorem{theorem}{Theorem}[section]
\newtheorem{lemma}[theorem]{Lemma}
\newtheorem{conj}[theorem]{Conjecture}
\newtheorem{prop}[theorem]{Proposition}
\newtheorem{corr}[theorem]{Corollary}

\makeatletter
\begingroup
\gdef\th@upshape{\normalfont
  \def\@begintheorem##1##2{%
        \item[\hskip\labelsep \theorem@headerfont ##1\ ##2.]}%
\def\@opargbegintheorem##1##2##3{%
   \item[\hskip\labelsep \theorem@headerfont ##1\ ##2\ (##3).]}}
\endgroup
\makeatother

\theoremstyle{upshape}

\newtheorem{defn}[theorem]{Definition}
\newtheorem{remark}[theorem]{Remark}
\newtheorem{exa}[theorem]{Example}

\makeatletter
\renewcommand{\subsection}{\@startsection{subsection}{2}{0pt}{-3ex
plus -1ex minus -0.2ex}{-2mm plus -0pt minus
-2pt}{\normalfont\bfseries}} 
\renewcommand{\subsubsection}{\@startsection{subsubsection}{3}{0pt}{-3ex
plus -1ex minus -0.2ex}{-2mm plus -0pt minus
-2pt}{\normalfont\bfseries}} 
\makeatother

\makeatletter
\@addtoreset{equation}{section}
\makeatother

\newcommand{\cntrct}                
{\hspace{2pt}\raisebox{1pt}{\text{$\lrcorner$}}\hspace{2pt}}

\newcommand{\proof}[1][Proof.]{\smallskip\noindent{\em #1}}
\def\endproof{\hfill\ensuremath{\square}\par\medskip}

\def\eqref#1{\thetag{\ref{#1}}}

\let\latexref=\ref
\def\ref#1{{\normalfont{\latexref{#1}}}}

\newcommand{\wt}{\widetilde}
\newcommand{\wh}{\widehat}

\newcommand{\idot}{{\:\raisebox{1pt}{\text{\circle*{1.5}}}}}
\newcommand{\hdot}{{\:\raisebox{3pt}{\text{\circle*{1.5}}}}}
\newcommand{\Z}{{\mathbb Z}}
\newcommand{\Q}{{\mathbb Q}}

\newcommand{\R}{{\mathbb R}}
\newcommand{\C}{{\mathbb C}}

\renewcommand{\phi}{\varphi}

\newcommand{\Fr}{{\sf Fr}}

\newcommand{\vH}{\check{H}}
\newcommand{\HH}{{\mathcal{H}}}

\def\dlim_#1{{\displaystyle\lim_{#1}}^\hdot}
\def\hocolim{\operatorname{\sf hocolim}}
\def\holim{\operatorname{\sf holim}}

\newcommand{\Hom}{\operatorname{Hom}}

\newcommand{\RHom}{\operatorname{RHom}}
\newcommand{\Tor}{\operatorname{Tor}}

\newcommand{\id}{\operatorname{\sf id}}

\newcommand{\gr}{\operatorname{\sf gr}}

\newcommand{\A}{{\cal A}}
\newcommand{\D}{{\cal D}}
\newcommand{\DM}{{\cal D}{\cal M}}

\newcommand{\B}{{\cal B}}

\newcommand{\Qq}{{\cal Q}}

\newcommand{\Sets}{\operatorname{Sets}}

\newcommand{\Maps}{\operatorname{Maps}}

\newcommand{\Iso}{{\operatorname{Iso}}}

\newcommand{\hhom}{\operatorname{\sf Hom}}
\newcommand{\Ind}{\operatorname{\sf Ind}}

\newcommand{\Sthom}{\operatorname{\sf StHom}}
\newcommand{\Stalg}{\operatorname{\sf StAlg}}

\newcommand{\St}{\operatorname{\sf St}}

\newcommand{\amod}{{\text{\rm -mod}}}

\newcommand{\ppt}{{\sf pt}}

\newcommand{\lotimes}{\overset{\sf\scriptscriptstyle L}{\otimes}}

\newcommand{\Spec}{\operatorname{Spec}}

\newcommand{\cchar}{\operatorname{\sf char}}

\newcommand{\Ab}{\operatorname{Ab}}

\newcommand{\M}{\operatorname{\mathcal M}}
\newcommand{\MM}{\operatorname{{\mathcal M}{\mathcal M}}}

\newcommand{\DML}{\D\!\M\!\Lambda}

\newcommand{\DLR}{\D\!\Lambda\text{{\upshape R}}}

\newcommand{\DF}{\D{\mathcal F}}

\newcommand{\FDM}{{\mathcal F}\D\!\M}
\newcommand{\DFDM}{\D{\mathcal F}\D\!\M}

\newcommand{\Real}{\operatorname{\sf Real}}
\newcommand{\Rep}{\operatorname{Rep}}
\newcommand{\Gal}{\operatorname{Gal}}

\renewcommand{\P}{{\mathbb P}}
\newcommand{\Ss}{{\mathbb S}}

\newcommand{\Dgalg}{\operatorname{\sf DG-Alg}}

\newcommand{\Dst}{\operatorname{{\mathcal D}^-{\sf st}}}
\newcommand{\Dgcomm}{\operatorname{\sf DG-Comm}}
\newcommand{\Comm}{\operatorname{\sf Comm}}
\newcommand{\Rspec}{\operatorname{R\Spec}}
\newcommand{\Mor}{\operatorname{\sf Morita}}
\newcommand{\Perf}{\operatorname{\sf Perf}}
\newcommand{\Alg}{\operatorname{\sf Alg}}
\newcommand{\Top}{\operatorname{\rm Top}}
\newcommand{\Spaces}{\operatorname{\rm Spaces}}

\newcommand{\EM}{\operatorname{\sf EM}}

\newcommand{\Symm}{\operatorname{{\sf Sym}}}

\newcommand{\THH}{\operatorname{THH}}
\newcommand{\TC}{\operatorname{TC}}

\newcommand{\I}{{\mathbb I}}
\newcommand{\can}{\operatorname{\sf can}}

\title{Motivic structures in non-commutative geometry}

\author{D. Kaledin}

\begin{document}

\maketitle

\tableofcontents

\section{Generalities on mixed motives.}\label{mm.sec}

The conjectural category $\MM$ of {\em mixed motives}, as described
by Deligne, Beilinson and others, unifies and connects various
cohomology theories which appear in modern algebraic
geometry. Recall that one expects $\MM$ to be a symmetric tensor
abelian category with a distiguished invertible object $\Z(1)$
called the {\em Tate motive}. One expects that for any smooth
projective algebraic variety $X$ defined over $\Q$, there exist
a functorial {\em motivic cohomology complex} $H^\hdot(X) \in
\D^b(\MM)$ with values in the derived category $\D^b(\MM)$, whose
cohomology groups
$$
H^i(X) \in \MM
$$
are called {\em motivic cohomology groups}. If $X$ is the projective
space $\P^n$, $n \geq 1$, then one expects to have
$$
H^{2i}(\P^n) \cong \Z(i)
$$
for $0 \leq i \leq n$, and $0$ otherwise. For a general $X$ and any
integer $j$, one defines the {\em absolute cohomology complex} by
$$
H^\hdot_{abs}(X,\Z(j)) = \RHom^\hdot_{\MM}(\Z(-j),H^\hdot(X)),
$$
with its cohomology groups $H^i_{abs}(X,\Z(j))$ known as {\em
absolute cohomology groups}. It is expected that the absolute
cohomology groups are related to the algebraic $K$-theory groups
$K_\idot(X)$ by means of a functorial {\em regulator map}
\begin{equation}\label{reg.abs}
r:K_\idot(X) \to \bigoplus_j H^{2j-\hdot}_{abs}(X,\Z(j)),
\end{equation}
and it is expected that the regulator map is ``not far from an
isomorphism'' (for example, it ought to be an isomorphism modulo
torsion).

\bigskip

The above picture, with its many refinements which we will not need,
is, unfortunately, still conjectural. In practice, one has to be
content with {\em categories of realizations}. These follow the same
general pattern, but the hypothetical category $\MM$ is replaced
with a known category $\Real$ whose definition axiomatizes the
features of a particular known cohomology theory. The prototype
example is that of $l$-adic cohomology. Recall that for any
algebraic variety $X/\Q$, its $l$-adic \'etale cohomology groups
$$
H_{et}^i(X,\Q_l)
$$
are $\Q_l$-vector spaces equipped with an additional structure of an
{\em $l$-adic representation} of the Galois group
$\Gal(\overline{\Q}/\Q)$. These representations form a tensor
symmetric abelian category $\Rep_l(\Gal(\overline{\Q}/\Q))$ with a
distiguished Tate module $\Q_l(1)$, and one can treat $l$-adic
cohomology as taking values in this category. One can them define a
double-graded absolute cohomology theory
$$
H^\hdot_{abs}(X,\Q_l(j)) =
H^\hdot(\Gal(\overline{\Q}/\Q),H^\hdot_{et}(X,\Q_l(j))),
$$
known as {\em absolute $l$-adic cohomology}, and construct a
regulator map of the form \eqref{reg.abs}. Conjecturally, we have an
exact tensor ``realization functor'' $\MM \to
\Rep_l(\Gal(\overline{\Q}/\Q))$, $l$-adic cohomology is obtained by
applying realization to motivic cohomology, and the \'etale
regulator map factors through the motivic one. In practice, one can
treat $\Rep_l(\Gal(\overline{\Q}/\Q))$ as a replacement for $\MM$,
and hope that the regulator map still captures essential information
about $K_\idot(X)$.

\medskip

In this paper, we will be concerned with another family of
cohomology theories and realizations which appear as refinements of
{\em de Rham cohomology}. By its very nature, de Rham cohomology of
a smooth algebraic variety $X$ has coefficients in the field or ring
of definition of $X$. Thus it is not necessary to require that $X$
is defined over $\Q$, and it is convenient to classify de Rham-type
cohomology theories by their rings of definitions. There are two
main examples.
\begin{enumerate}
\item The ring of definition is either $\R$ or $\C$; the
  corresponding category of realizations is Deligne's category of
  mixed $\R$-Hodge structures, and the absolute cohomology theory is
  Hodge-Deligne cohomology (with a refinement by Beilinson). The
  regulator map is the subject of famous Beilinson Conjectures.
\item The ring of definition is $\Z_p$; the corresponding category
  of realizations is the category of {\em filtered Dieudonn\'e
  modules} of Fontaine-Lafaille \cite{FL}, and the absolute
  cohomology theory is {\em syntomic cohomology} of Fontaine and
  Messing \cite{FM}.
\end{enumerate}
The goal of this paper is to report on recent discoveries and
conjectures which state, roughly speaking, that all these additional
``motivic'' structures on de Rham cohomology of an algebraic variety
should exist in a much more general setting of periodic cyclic
homology of properly understood {\em non-commutative} algebraic
varieties. As opposed to the usual commutative setting, the
``classical'' case \thetag{i} is more difficult and largely
conjectural; in the $p$-adic case \thetag{ii}, most of the
statements have been proved. Moreover, the $p$-adic story shows an
unexpected relation to algebraic topology which we will also
explain. Before we start, however, we should define exactly what we
mean by a ``non-commutative algebraic variety'', and recall basic
facts on cyclic homology.

\section{Non-commutative setting.}\label{hc.sec}

We start by a brief recollection on cyclic homology; a very good
overview can be found in J.-L. Loday's book \cite{Lo}, and an old
overview \cite{FT} is also quite useful. {\em Hochschild homology}
$HH_\idot(A/k)$ of an associative unital algebra $A$ flat over a
commutative ring $k$ is given by
$$
HH_\idot(A) = HH_\idot(A/k) = \Tor_\idot^{A^{opp} \otimes_k A}(A,A),
$$
where $A^{opp}$ is $A$ with multiplication written in the opposite
direction. It has been discovered by Hochschild, Kostant and
Rosenberg \cite{HKR} that if $A$ is commutative and $X = \Spec A$ is
a smooth algebraic variety over $k$, then
$$
HH_i(A) \cong H^0(X,\Omega^i(X)),
$$
the space of $i$-forms on $X$ over $k$. {\em Cyclic homology}
$HC_\idot(A)$ is a refinement of Hochschild homology discovered
independently by A. Connes and B. Tsygan. It is functorial in $A$,
and related to $HH_\idot(A)$ by the {\em Connes' long exact
sequence}
$$
\begin{CD}
HH_\idot(A) @>>> HC_\idot(A) @>{u}>> HC_{\idot-2}(A) @>>>,
\end{CD}
$$
where $u$ is a canonical {\em periodicity map} of degree $2$. Both
$HH_\idot(A)$ and $HC_\idot(A)$ can be represented by functorial
complexes $CH_\idot(A)$, $CC_\idot(A)$, and the Connes' exact
sequence then becomes a short exact sequence of complexes.  The
complex $CC_\idot(A)$ is the total complex of a bicomplex
\begin{equation}\label{hc.eq}
\begin{CD}
\dots @>>> A @>{\id}>> A @>{0}>> A\\
@. @AA{b}A @AA{b'}A @AA{b}A \\
\dots @>>> A \otimes A @>{\id + \tau}>> A \otimes A 
@>{\id - \tau}>> A \otimes A\\
@. @AA{b}A @AA{b'}A @AA{b}A\\
\dots @. \dots @. \dots @. \dots\\
@. @AA{b}A @AA{b'}A @AA{b}A \\
\dots @>>> A^{\otimes n} @>{\id + \tau + \dots + \tau^{n-1}}>> 
A^{\otimes n} @>{\id - \tau}>> A^{\otimes n}\\
@. @AA{b}A @AA{b'}A @AA{b}A
\end{CD}
\end{equation}
Here it is understood that the whole thing extends indefinitely to
the left, all the even-numbered columns are the same, all
odd-numbered columns are the same, and the bicomplex is invariant
with respect to the horizontal shift by $2$ columns which gives the
periodicity map $u$. The map $\tau:A^{\otimes i} \to \A^{\otimes i}$
is the cyclic permutation of order $i$ multiplied by $(-1)^{i+1}$,
and $b$, $b'$ are certain explicit differentials expressed in terms
of the multiplication map $m:A^{\otimes 2} \to A$.  The complex
$CH_\idot(A)$ is the rightmost column of \eqref{hc.eq}, and also any
odd-numbered column when counting from the right; the even-numbered
columns are acyclic.

{\em Periodic cyclic homology} $HP_\idot(A)$ is obtained by
inverting the map $u$, namely, $HP_\idot(A)$ is the homology of the
complex
$$
CP_\idot(A) = \lim_{\overset{u}{\gets}}CC_\idot(A)
$$
(explictily, this is the total complex of a bicomplex obtained by
extending \eqref{hc.eq} to the right as well as to the left). {\em
Negative cyclic homology} $HC^-_\idot(A)$ is the homology of the
complex $CC^-_\idot(A)$ obtained as the third term in a short exact
sequence
$$
\begin{CD}
0 @>>> CC^-(A) @>>> CP_\idot(A) @>>> CC_{\idot-2}(A) @>>> 0
\end{CD}
$$
(equivalently, one extends \eqref{hc.eq} to the right but not to the
left).

\medskip

The reason cyclic homology is interesting in algebraic geometry is
the following comparison theorem. In the situation of the
Hochschild-Kostant-Rosenberg Theorem, let $d$ be the dimension of $X
= \Spec A$, and assume in addition that $d!$ is invertible in the
base ring $k$. Then there exists a canonical isomorphism
\begin{equation}\label{hp.dr}
HP_\idot(A) \cong H^\hdot_{DR}(X)((u)),
\end{equation}
where the right-hand side is a shorthand for ``formal Laurent power
series in one variable $u$ of degree $2$ with coefficients in de
Rham cohomology $H_{DR}(X)$''.

By \eqref{hp.dr}, periodic cyclic homology classes can be thought of
as non-commutative generalizations of de Rham cohomology
classes. Some information is lost in this generalization: because of
the presense of $u$ in the right-hand side of \eqref{hp.dr}, what we
recover from $HP_\idot(A)$ is not the de Rham cohomology of $X$ but
rather, the de Rham cohomology of the product $X \times \P^\infty$
of $X$ and the infinite projective space $\P^\infty$, where we
moreover invert the generator $u \in H^2_{DR}(\P^\infty)$. Thus
given a category of realizations $\Real$ and a $\Real$-valued
refinement of de Rham cohomology, the appropriate target for its
non-commutative generalization is not the derived category
$\D(\Real)$ but the {\em twisted $2$-periodic derived category}
$\D^{per}(\Real)$ obtained by inverting quasiisomorphisms in the
category of complexes $M_\idot$ of objects in $\Real$ equipped with
an isomorphism $u:M_\idot \cong M_\idot(-1)[2]$, where we denote
$M(n) = M \otimes \Z(n)$, $n \in \Z$.

We note, however, that this causes no problem with the regulator
map, since the summation in the right-hand side of \eqref{reg.abs}
is the same as in the right-hand side of \eqref{hp.dr}. Thus for a
$\Real$-valued refinement $H^\hdot_{\Real}(-)$ of de Rham cohomology
and any smooth affine algebraic variety $X = \Spec A$, the regulator
map \eqref{reg.abs} takes the form
$$
K_\idot(A) \to \RHom^\hdot_{\D^{per}(\Real)}(k,HP_\idot(A)) =
\RHom^\hdot_{\D^{per}(\Real)}(k,H^\hdot_{\Real}(X)((u))),
$$
where $k$ in the right-hand side is the unit object of $\Real$.

\medskip

Somewhat surprisingly, non-affine algebraic varieties can be
included in the above picture with very little additional effort. To
do it, it is convenient to use the machinery of differential graded
(DG) algebras and DG categories. An excellent overview can be found
in \cite{kel}; for the convenience of the reader, let us summarize
the relevant points.

Roughly speaking, a $k$-linear DG category is a category $C^\hdot$
whose $\Hom$-sets $C^\hdot(-,-)$ are equipped with a structure of
complexes of $k$-modules in such a way that composition maps are
$k$-linear and compatible with the differentials (for precise
definitions, see \cite[Section 2]{kel}). For any small $k$-linear DG
category $C^\hdot$, one defines a triangulated {\em derived category
  of DG modules} $\D(C^\hdot)$ (\cite[Section 3]{kel}). Any
$k$-linear DG functor $\gamma:C_1^\hdot \to C_2^\hdot$ induces a
trinagulated functor $\gamma^*:\D(C_2^\hdot) \to \D(C_1^\hdot)$. The
functor $\gamma$ ia a {\em derived Morita equivalence} if the
induced functor $\gamma^*$ is an equivalence of triangulated
categories. It turns out -- this mostly due to the work of
G. Tabuada and B. To\"en, see \cite[Section 4]{kel} and references
therein -- that there is a closed model structure on the category of
small $k$-linear DG categories whose weak equivalences are exactly
derived Morita equivalences. Denote by $\Mor(k)$ the corresponding
homotopy category, that is, the category of ``small $k$-linear DG
categories up to a derived Morita equivalence''.

Any $k$-algebra $A$ is a $k$-linear DG category with one object
$\ppt$ and $\Hom(\ppt,\ppt)=A$ placed in degree $0$, so that we have
an embedding $\Alg(k) \to \Mor(k)$ from the category $\Alg(k)$ of
associative $k$-algebras to $\Mor(k)$. Then, as explained in
\cite[Section 5]{kel}, Hochschild homology, cyclic homology,
periodic cyclic homology and negative cyclic homology extend to
functors
$$
\Mor(k) \to \D(k).
$$
Moreover, so does the algebraic $K$-theory functor $K^\hdot(-)$, and
other ``additive invariants'' in the sense of \cite[Section 5]{kel}.

In general, a DG category with one object $\ppt$ is the same thing
as an associative unital DG algebra $A^\hdot =
\Hom^\hdot(\ppt,\ppt)$. The category of DG algebras over $k$ has a
natural closed model structure whose weak equivalences are
quasiisomorphisms, and whose fibrations are surjective maps. The
corresponding homotopy category $\Dgalg(k)$ is the category of DG
algebras ``up to a quasiisomorphism''. One shows that a
quasiisomorphism between DG algebras is in particular a derived
Morita equivalence, so that we have a natural functor
\begin{equation}\label{dg.mor}
\Dgalg(k) \to \Mor(k).
\end{equation}
It is not diffucult to show that for every cofibrant DG algebra
$A^\hdot$, the individual terms of the complex $A^\hdot$ are flat
$k$-modules. In this case, the Hochschild, cyclic etc. homology of
$A^\hdot$ are especially simple -- they are given by exactly the
same bicomplex \eqref{hc.eq} and its versions as in the case of
ordinary algebras. This is manifestly invariant
under quasiisomorphisms, so that the Hochschild, cyclic etc. homology
obviously descend to functors from $\Dgalg(k)$ to the derived
category $\D(k)$. The DG category approach shows that there is even
more invariance: even if two DG algebras $A_1^\hdot$, $A_2^\hdot$
are not quasiisomorphic but only have isomorphic images in
$\Mor(k)$, their Hochschild, cyclic etc. homology is naturally
identified. This statement is already non-trivial in the case of
usual algebras, see \cite[Section 1.2]{Lo}.

\begin{defn}
  A DG category $T_\idot \in \Mor(k)$ is {\em derived-affine} if it
  lies in the essential image of the functor \eqref{dg.mor}.
\end{defn}

\begin{remark}
  A small $k$-linear DG category $C^\hdot$ with a finite number of
  objects is automatically derived-Morita equivalent to a DG algebra
  $A^\hdot$, thus affine. For example, one can take
$$
A^\hdot = \bigoplus_{c,c'}C^\hdot(c,c'),
$$ 
where the sum is taken over all pairs of objects in $\C^\hdot$.
\end{remark}

Now, it has been proved (\cite{BV} combined with \cite{kelold}) that
for any quasiseparated quasicompact scheme $X$ over $k$, there
exists a DG algebra $A^\hdot/k$ such that the derived category
$\D(X)$ of quasicoherent sheaves on $X$ is equivalent to the derived
category $\D(A^\hdot)$,
$$
\D(X) \cong \D(A^\hdot),
$$
and such a DG algebra $A^\hdot$ is unique up to a derived Morita
equivalence, so that we have a canonical functor from the category
of algebraic varieties over $k$ to the category $\Mor(k)$. Roughly
speaking, any algebraic variety is derived Morita-equivalent to a DG
algebra, or, in a succint formulation of \cite{BV}, ``every
algebraic variety is derived-affine''.

Moreover, it turns out that the properties of $X$ which are relevant
for the present paper are reflected in the properties of a
Morita-equivalent DG algebra $A^\hdot$. For example, one introduces
the following (see e.g. \cite{KS}).

\begin{defn}
\begin{enumerate}
\item A DG algebra $A^\hdot/k$ is {\em proper} if $A^\hdot$ is
  perfect as an object in the derived category $\D(k)$ of complexes
  of $k$-modules.
\item A DG algebra $A^\hdot/k$ is {\em smooth} if $A^\hdot$ is
  perfect as an object in the derived category $\D(A^{\hdot
  opp}\otimes A^\hdot)$ of $A^\hdot$-bimodules.
\end{enumerate}
\end{defn}

Then $A^\hdot$ is proper, resp. smooth if and only if $X$ is proper,
resp. smooth (in the affine case $X = \Spec A$, the second claim is
the famous Serre regularity criterion). Moreover, the correspondence
$X \mapsto A^\hdot$ is compatible with algebraic $K$-theory,
$K_\idot(X) \cong K_\idot(A^\hdot)$, and if the variety $X/k$ is
smooth of dimension $d$, and $d!$ is invertible in $k$, then the
Hochschild homology of such a Morita-equivalent DG algebra $A^\hdot$
is canonically isomorphic to
$$
HH_i(A^\hdot) \cong \bigoplus_j H^j(X,\Omega^{i+j}(X)),
$$
the so-called ``Hodge cohomology'' of $X$, while the periodic cyclic
homology $HP_\idot(A)$ is exactly as in \eqref{hp.dr}. 

Thus as far as homological invariants are concerned, one can treat
DG algebras ``up to a derived Morita-equivalence'' as
non-commutative generalizations of algebraic varieties:
\begin{itemize}
\item A non-commutative algebraic variety over $k$ is a DG algebra
  $A^\hdot$ over $k$ considered as an object of the Tabuada-To\"en
  category $\Mor(k)$.
\end{itemize}
This is the point of view we will adopt.

\section{Hodge-to-de Rham spectral sequence.}\label{hdr.sec}

A convenient way to pack all the structures related to Hochschild
homology $HH_\idot(A^\hdot)$ of a DG algebra $A^\hdot/k$ is by
considering the equivariant derived category $\D_{S^1}(k)$ of
$S^1$-equivariant constructible sheaves of $k$-modules on the point
$\ppt$. Then the claim is that the Hochschild homology complex
$CH_\idot(A^\hdot)$, while {\em a priori} simply a complex of
$k$-modules, in fact underlies a canonical object
$\wt{CH}_\idot(A^\hdot) \in \D_{S^1}(k)$ (loosely speaking,
``$CH_\idot(A^\hdot)$ carries a canonical $S^1$-action'').  The
negative cyclic homology appears as $S^1$-equivariant cohomology
$$
H^\hdot_{S^1}(\ppt,\wt{CH}_\idot(A^\hdot)),
$$
the periodicity map $u$ is the generator of $H^\hdot_{S^1}(\ppt)
\cong H^\hdot(BS^1)$, and $HP_\idot(A^\hdot)$ is the localization
$HC_\idot^-(A^\hdot)(u^{-1})$.

Another way to pack the same data is by considering the {\em
  filtered derived category} $\DF(k)$ of $k$-modules of \cite{BBD}
-- that is, the triangulated category obtained by considering
complexes $V_\idot$ of $k$-modules equipped with a descreasing
filtration $F^\hdot$ numbered by all integers, and inverting those
maps which induce quasiisomorphisms on the associated graded
quotients $\gr^F$. This has a ``twisted $2$-periodic'' version
$\DF^{per}(k)$, obtained from filtered complexes $V_\idot$ equipped
with an isomorphism $V_\idot \cong V_\idot[2](-1)$, where $(-1)$
means renumbering the filtration: $F^iV(-1) = F^{\hdot+1}V$.

\begin{lemma}\label{trivial.lemma} We have
$$
\D_{S^1}(k) \cong \DF^{per}(k).
$$
\end{lemma}

\proof[Sketch of the proof.] Let us just indicate the equivalence:
it sends $V_\idot \in \D_{S^1}(k)$ to the equivariant cohomology
complex $C^\hdot_{S^1}(\ppt,V_\idot)(u^{-1})$, with the
(generalized) filtration given by
$$
F^iH^\hdot_{S^1}(\ppt,V_\idot)(u^{-1}) =
u^iC^\hdot_{S^1}(\ppt,V_\idot),
$$
where $u \in C^2_{S^1}(k)$ represents the generator of the
equivariant cohomology ring $H^\hdot_{S^1}(\ppt,k) \cong k[u]$.
\endproof

In the case of the Hochschild homology complex
$\wt{CH}_\idot(A^\hdot)$, the corresponding periodic filtered
complex is $CP_\idot(A^\hdot)$, with the filtration given by
$$
F^iCP_\idot((A^\hdot) = u^iCC^-_\idot(A^\hdot) \subset
CP_\idot(A^\hdot).
$$

\medskip

One can treat $\DF^{per}(k)$ as a very crude ``category of
realization'' $\Real$ in the sense of Section~\ref{mm.sec}, or
rather, of its periodic derived category $\D^{per}(\Real)$. The
expected regulator map then takes the form
\begin{equation}\label{non-comm.reg}
K_\idot(A^\hdot) \to HC^-_\idot(A^\hdot) =
\RHom^\hdot_{\DF^{per}(k)}(k,HP_\idot(A)).
\end{equation}
Such a map does indeed exist, see \cite[Chapter 8]{Lo}. In general,
it is very far from being an isomorphism. The only general result is
a theorem of T. Goodwillie \cite{good} which shows that at least the
tangent spaces to both sides are the same. Namely, given an alegbra
$A$ with an ideal $I \subset A$, one defines the relative $K$-theory
$K_\idot(A,I)$ spectrum as the cone of the natural map $K_\idot(A)
\to K_\idot(A/I)$, and analogously for the cyclic homology
functors. Then it has been proved in \cite{good} that if $k$ is a
field of characteristic $0$ and $I \subset A$ is a nilpotent ideal,
then the map
$$
K_\idot(A,I) \to HC^-(A,I)
$$
induced by the regulator map \eqref{non-comm.reg} is a
quasiisomorphism. An analogous statement also holds for DG algebras
over $k$.

\medskip

While filtered complexes are a very crude approximation to mixed
motives, already on this level the smoothness and properness of a DG
algebra leads to non-trivial consequences. Namely, a filtered
complex gives rise to a spectral sequence. In the case of cyclic
homology, it takes the form
\begin{equation}\label{h.d.r}
HH_\idot(A^\hdot)((u)) \Rightarrow HP_\idot(A^\hdot),
\end{equation}
where we use the same shorthand as in \eqref{hp.dr}. When the DG
algebra $A^\hdot/k$ is Morita-equivalent to a smooth algebraic
variety $X/k$, the filtration $F^\hdot$ on $HP_\idot(A^\hdot) \cong
H^\hdot_{DR}(X)((u))$ is just the Hodge filtration on de Rham
cohomology, and \eqref{h.d.r} is the usual Hodge-to-de Rham spectral
sequence
$$
H^p(X,\Omega^q(X))\Rightarrow H^{p+q}_{DR}(X)
$$
tensored with $k((u))$. Because of this, \eqref{h.d.r} in general is
also called ``Hodge-to-de Rham spectral sequence''. Then the
following is a partial proof of a general conjecture of
M. Kontsevich and Ya. Soibelman \cite{KS1}.

\begin{theorem}[\cite{K}]\label{hdr.thm}
Assume that $A^\hdot$ is a smooth and proper DG algebra over a field
$k$ of characteristic $\cchar k = 0$. Assume further that $A^i = 0$
for $i < 0$. Then the Hodge-to-de Rham spectral sequence
\eqref{h.d.r} degenerates.
\end{theorem}

The assumption $A^i = 0$, $i < 0$ is technical (note, however, that
it can always be achieved for a DG algebra $A^\hdot$ corresponding
to a smooth and proper algebraic variety $X/k$, see
e.g. \cite[Theorem 4]{O}).

In the usual commutative case, the Hodge-to-de Rham degeneration
statement is well-known and has two proofs. Classically, it follows
from the general complex-analytic package of Hodge theory and
harmonic forms. An alternative proof by Deligne and Illusie
\cite{DI} uses reduction to positive characteristic and $p$-adic
methods. So far, it is only the second technique that has been
generalized to the non-commutative case. We will now explain this.

\section{Review of Filtered Dieudonn\'e modules.}\label{fdm.1.sec}

A $p$-adic analog of the notion of a mixed Hodge structure has been
introduced in 1982 by Fontaine and Lafaille \cite{FL}. Here is the
definition.

\begin{defn}\label{fdm.defn}
  Let $k$ be a finite field of characteristic $p$, with its
  Frobenius map, and let $W$ be its ring of Witt vectors, with its
  canonical lifting $\phi$ of the Frobenius map. A {\em filtered
    Dieudonn\'e module} over $W$ is a finitely generated $W$-module
  $M$ equipped with a decreasing filtration $F^\hdot M$, indexed by
  all integers and such that $\cap F^iM=0$, $\cup F^iM=M$, and a
  collection of Frobenius-semilinear maps $\phi_i:F^iM \to M$, one
  for each integer $i$, such that
\begin{enumerate}
\item $\phi_i|_{F^{i+1}M} = p \phi^{i+1}$, and
\item the map
$$
\sum \phi_i:\bigoplus_iF^iM \to M
$$
is surjective.
\end{enumerate}
\end{defn}

We will denote by $\FDM(W)$ the category of filtered Dieudonn\'e
modules over $W$. It is an abelian category. A symmetric tensor
product in $\FDM(W)$ is defined in the obvious way, and we have the
Tate object $W(1)$ given by: $W(1) = W$ as a $W$-module, $F^1W(1) =
W(1)$, $F^2W(1)=0$, $\phi_1:F^1W(1) \to W(1)$ equal to $\phi$. We
also have the derived category $\D(\FDM(W))$.

If a filtered Dieudonn\'e module $M \in \FDM(W)$ is annihilated by
$p$, then \thetag{i} of Definition~\ref{fdm.defn} insures that the
map in \thetag{ii} factors through a surjective map
$$
\wt{\phi}:\gr_F^\hdot M \to M.
$$
Since both sides are $k$-vector spaces of the same dimension,
$\wt{\phi}$ must be an isomorphism. For a general filtered
$W$-module $\langle M,F^\hdot \rangle$, one lets $\wt{M}$ be the
cokernel of the map
\begin{equation}\label{wt.M}
\begin{CD}
\bigoplus_i F^iM @>{t - p\id}>> \bigoplus_i F^iM,
\end{CD}
\end{equation}
where $t:F^{\hdot+1}M \to F^\hdot M$ is the tautological
embedding. Then again, \thetag{i} insures that the map
$\sum_i\phi_i$ factors through a map
\begin{equation}\label{wt.phi}
\wt{\phi}:\wt{M} \to M 
\end{equation}
and this map must be an isomorphism if \thetag{ii} were to be
satisfied. This allows to generalize the definition of a filtered
Dieudonn\'e module: instead of a finitely generated filtered
$W$-module, one can consider a filtered $W$-module $\langle
M,F^\hdot \rangle$ such that $M$ is $p$-adically complete and
complete with respect to the topology induced by $F^\hdot$ (these
conditions together with the non-degeneracy conditions $\cap
F^iM=0$, $\cup F^iM=M$ insure that the map \eqref{wt.M} is
injective). Then a {\em unbounded Dieudonn\'e module} structure on
$M$ is given by a Frobenius-semilinear isomorphism $\wt{\phi}$ of
the form \eqref{wt.phi}.

I do not know whether the category of unbounded filtered Dieudonn\'e
modules is still abelian. However, complexes of unbounded filtered
Dieudonn\'e modules can be defined in the obvious way, and the
correspondence $M \mapsto \wt{M}$ sends filtered quasiisomorphisms
into quasiisomorphisms, so that we obtain a triangulated derived
category $\DFDM(W) \supset \D(\FDM(W))$ and its twisted $2$-periodic
version $\DFDM^{per}$.

Moreover, one can drop the requirement that the map $\wt{\phi}$ is
an isomorphism and allow it to be an arbitrary map. Let us call the
resulting objects ``weak filtered Dieudonn\'e modules''. The
category of weak filtered Dieudonn\'e modules is definitely not
abelian, but the above procedure stil applies: we can invert
filtered quasiisomorphisms and obtain triangulated categories
$\wt{\DFDM}(W)$, $\wt{\DFDM}^{per}(W)$. We then have a fully
faithful inclusions $\DFDM(W) \subset \wt{\DFDM}(W)$,
$\DFDM^{per}(W) \subset \wt{\DFDM}^{per}(W)$, and one can show that
their essential images consist of those $M_\idot$ in
$\wt{\DFDM}(W)$, resp. $\wt{\DFDM}^{per}(W)$ for which the map
$\wt{\phi}$ of \eqref{wt.phi} is a quasiisomorphism.

\medskip

Assume given a algebraic variety $X$ smooth over $W$, of dimension
$d < p$. Then de Rham cohomology $H^\hdot_{DR}(X/W)$ equipped with
the filtration induced by the stupid filtration on the de Rham
complex has the structure of a complex of generalized filtered
Dieudonn\'e modules. If $X/W$ is proper, the groups $H^i_{DR}(X/W)$
are finitely generated, so that they are filtered Dieudonn\'e
modules in the strict sense (and the filtration is then the Hodge
filtration). This Dieudonn\'e module structure can be seen
explicitly under the following strong additional assumption:
\begin{itemize}
\item the Frobenius endomorphism $\Fr$ of the special fiber $X_k = X
  \otimes_W k$ of $X/W$ lifts to a Frobenius-semilinear endomorphism
  $\wt{\Fr}:X \to X$.
\end{itemize}
Then one checks easily that for any $i \geq 0$, the natural map
$\wt{\Fr}^*:\Omega^i(X/W) \to \Omega^i(X/W)$ is divisible by
$p^i$. The Dieudonn\'e module structure maps $\phi_i$ are induced by
the corresponding maps $\frac{1}{p^i}\wt{\Fr}^*$. We note that in
this special case, the map $\phi_i$ sends $F^i$ into $F^i$. In the
general case, the construction is due to G. Faltings \cite[Theorem
4.1]{falt}; roughly speaking, it uses a comparison theorem which
gives a quasiisomorphism
$$
H^\hdot_{cris}(X_k) \cong H^\hdot_{DR}(X),
$$
where in the left-hand side, we have the cristalline cohomology of
the special fiber $X_k$. The Frobenius endomorphism of $X_k$ induces
an endomorphism on cristalline cohomology, and this gives the
structure map $\phi_0$. By an additional argument, one shows that
$\phi_0|F^i$ is canonically divisible by $p^i$, and this gives the
other structure maps $\phi_i$ (in general, they do not preserve the
Hodge filtration $F^\hdot$).

In particular, for any smooth $X/W$, one has the isomorphism
\eqref{wt.phi}. Its reduction $\mod p$ is an isomorphism
\begin{equation}\label{cartier}
\gr_F^\hdot H^\hdot_{DR}(X_k) \cong \bigoplus_i
H^{\hdot-i}(X_k,\Omega^i(X_k)) \cong H^\hdot_{DR}(X_k)
\end{equation}
between Hodge and de Rham cohomology of the special fiber $X_k$. If
$X$ is affine, this is nothing but the inverse to the Cartier
isomorphism, discovered by P. Cartier back in the 1950-ies; as such,
it depends only on the special fiber $X_k$ and not on the lifting
$X/W$. In the general case, it has been shown by Deligne and Illusie
in \cite{DI} that \eqref{cartier} depends on the lifting $X
\otimes_W W_2(k)$ of $X_k$ to the second Witt vectors ring $W_2(k) =
W(k)/p^2$ (but not on the lifting to higher orders, nor even on the
existence of such a lifting).

The absolute cohomology theory associated to the $\FDM$-valued
refinement of de Rham cohomology is the {\em syntomic cohomology} of
Fontaine and Messing. As it happens, the functors
$\RHom^\hdot(W(-j),-)$ in the category $\DFDM$ are easy to compute
explicitly --- for any complex $M_\idot \in \DFDM$,
$\RHom^\hdot(W,-)$ is the cone of the natural map
$$
\begin{CD}
F^jM_\idot @>{\id - \phi_j}>> M_\idot.
\end{CD}
$$
When applied to a smooth proper variety $X/W$, this gives syntomic
cohomology groups $H^\hdot_{synt}(X,\Z_p(j))$. The construction can
even be localized with respect to the Zariski topology on $X_k$, so
that the syntomic cohomology is expressed as hypercohomology of
$X_k$ with coefficients in certain canonical complexes of Zariski
sheaves, as in \cite{FM}.

The existence and properties of the regulator map for the syntomic
cohomology have been studied by M. Gros \cite{gros1, gros2}. In
principle, one can construct the regulator by the standard procedure
for ``twisted cohomology theories'' in the sense of \cite{BlO}, but
there is one serious problem: the filtered Dieudonn\'e module
structure on $H^\hdot_{DR}(X)$ only exists if $p > \dim X$. Since
the standard procedure works by considering infinite projective
spaces and Grassmann varieties, this condition is inevitably broken
no matter what $p$ we start with. To circumvent this, Gros had to
modify (in \cite{gros2}) the definition of syntomic cohomology by
including additional structures such as the rigid analytic space
associated to $X/W$. The resulting picture becomes extremely
complex, and at present, it is not clear whether it can be
generalized to non-commutative varieties.

\section{FDM in the non-commutative case.}\label{fdm.2.sec}

What we do have for non-commutative varieties is the following
result.

\begin{defn}
The {\em Hochschild cohomology} $HH^\hdot(A^\hdot/R)$ of a DG
algebra $A^\hdot$ over a ring $R$ is given by
$$
HH^\hdot(A^\hdot/R) = \RHom^\hdot_{A^{\hdot opp} \otimes_R
  A}(A^\hdot,A^\hdot).
$$
\end{defn}

\begin{theorem}[\cite{K}]\label{car}
Assume given an associative DG algebra $A^\hdot$ over a finite field
$k$. Assume that $A^i = 0$ for $i < 0$. Assume also that $A^\hdot$
is smooth, that it can be lifted to a flat DG algebra $\wt{A}^\hdot$
over $W_2(k)$, and that $HH^i(A^\hdot) = 0$ for $i \geq 2p-1$. Then
there exists a canonical Cartier-type isomorphism
$$
HH_\idot(A^\hdot)((u)) \cong HP_\idot(A^\hdot).
$$
\end{theorem}

\begin{remark} If a DG algebra $A^\hdot$ is derived
  Morita-equivalent to a smooth algebraic variety $X/k$, then we
  have $HH^i(A^\hdot) = 0$ automatically for $i > 2 \dim X$, so
  that the last condition on $A^\hdot$ in Theorem~\ref{car} reduces
  to the condition $p > \dim X$ already mentioned in
  Section~\ref{fdm.1.sec}.
\end{remark}

\begin{remark}
Theorem~\ref{hdr.thm} easily follows from Theorem~\ref{car} by the
same dimension argument as in the original proof of Deligne and
Illusie in \cite{DI}. The only non-trivial additional input is a
beautiful recent theorem of B. To\"en \cite{To} which claims that a
smooth and proper DG algebra $A^\hdot$ over a field $K$ comes from a
smooth and proper DG algebra $A^\hdot_R$ over a finitely generated
subring $R \subset K$, $A^\hdot \cong A^\hdot_R \otimes_R K$. This
allows one to reduce problems from $\cchar 0$ to $\cchar p$.
\end{remark}

Let us give a very rough sketch of how Theorem~\ref{car} is proved
(for more details, see \cite{goet}, and the complete proof in a
slightly different language is in \cite{K}). As in the commutative
story, there are two cases for Theorem~\ref{car}: the easy case when
one can construct the Cartier map explicitly, and the general
case. The easy case is when $A^\hdot = A$ is concentrated in degree
$0$, and the algebra $A$ admits a so-called {\em quasi-Frobenius
map}.

\begin{lemma}[\cite{K}]\label{tt.le}
For any vector space $V$ over the finite field $k$ of characteristic
$\cchar k = p > 0$, there is a canonical Frobenius-semilinear
isomorphism
$$
\vH^\hdot(\Z/p\Z,V) \cong \vH^\hdot(\Z/p\Z,V^{\otimes p}),
$$
where $\vH^\hdot(\Z/p\Z,-)$ means the Tate (co)homology of the group
$\Z/p\Z$, the action of $\Z/p\Z$ on $V$ is trivial, and the action
on $V^{\otimes p}$ is by the longest cycle permutation
$\sigma:V^{\otimes p} \to V^{\otimes p}$.\endproof
\end{lemma}

\begin{defn}[\cite{K}]
A {\em quasi-Frobenius map} for an algebra $A/k$ is a
$\Z/p\Z$-equivariant algebra map
$$
\Phi:A \to A^{\otimes p}
$$
which induces the standard isomorphism of Lemma~\ref{tt.le} on Tate
cohomology $\vH^\hdot(\Z/p\Z,-)$.
\end{defn}

If the algebra $A$ admits a quasi-Frobenius map $\Phi$, then the
construction of the Cartier isomorphism proceeds as follows. First,
recall that for any algebra $B$ equipped with an action of a group
$G$, the {\em smash product algebra} $B \# G$ is the group algebra
$B[G]$ but with the twisted product given by
$$
(b_1 \cdot g_1)(b_2 \cdot g_2) = b_1b_2^{g_1} \cdot g_1g_2,
$$
and one has a canonical decomposition
\begin{equation}\label{burgh}
HP_\idot(B \# G) = \bigoplus_{\langle g \rangle} HP_\idot(B \# G)_g
\end{equation}
into components numbered by conjugacy classes of elements in $G$
(these components are sometimes called {\em twisted sectors}). Next,
let $G$ be the cyclic group $\Z/p\Z$, and let $\sigma \in G$ be the
generator. Then one can show that if the $G$-action on $B$ is
trivial, then
\begin{equation}\label{tw.1}
HP_\idot(B \# G)_\sigma \cong \wt{HP_\idot(B)},
\end{equation}
where $HP_\idot(B)$ in the right-hand side is equipped with the
Hodge filtration, and $\wt{M}$ for a filtered group $M$ means the
cokernel of the map \eqref{wt.M}, as in Section~\ref{fdm.1.sec}. One
the other hand, if we take the $p$-th power $B^{\otimes p}$ with
$\sigma$ acting by the longest cycle permutation, then one can show
that 
\begin{equation}\label{tw.2}
HP_\idot(B^{\otimes p} \# G)_\sigma \cong HP_\idot(B).
\end{equation}
Both the isomorphisms \eqref{tw.1} and \eqref{tw.2} are completely
general and valid for algebras over any ring. So is the
decomposition \eqref{burgh}, which is moreover functorial with
respect to $G$-equivariant maps. We now apply this to our algebras
$A$ and $A^{\otimes p}$ over $k$, with the $G$-action as in
Lemma~\ref{tt.le}. The quasi-Frobenius map $\Phi$ induces a map
$$
\phi:\wt{HP_\idot(A)} \cong HP_\idot(A \# G)_\sigma \to
HP_\idot(A^{\otimes p} \# G)_\sigma \cong HP_\idot(A),
$$
and since $p$ annihilates $HP_\idot(A)$, we have
$$
\wt{HP_\idot(A)} \cong \gr_F^\hdot HP_\idot(A) \cong
HH_\idot(A)((u)).
$$
The map $\phi$ is the Cartier map of Theorem~\ref{car}. One then
shows that it is an isomorphism; this requires one to assume that
$A$ is smooth.

\medskip

The general case of Theorem~\ref{car} is handled by finding a
replacement for a quasi-Frobenius map; as far as the cyclic homology
is concerned, the argument stays the same. One first shows that for
any unital associative algebra $A/k$, there exists a completely
canonical diagram
$$
\begin{CD}
A @<{\alpha}<< Q_\idot(A) @>{\Phi}>> P_\idot(A) @<{\beta}<<
A^{\otimes p}
\end{CD}
$$
of DG algebras equipped with an action of $G=\Z/p\Z$ and
$G$-equivariant maps between them. The G action on $A$ and
$Q_\idot(A)$ is trivial. In addition, if $A$ is smooth, the map
$$
HP_\idot(A^{\otimes p} \# G)_\sigma \to HP_\idot(P_\idot(A) \#
G)_\sigma
$$
induced by the map $\beta$ is an isomorphism (although in general,
this isomorphism does not preserve the Hodge filtration). Thus as
before, $\Phi$ induces a canonical map
$$
\overline{\phi}:HH_\idot(Q_\idot(A))((u)) \cong
\wt{HP_\idot(Q_\idot(A))} \to HP_\idot(A).
$$
To construct the Cartier map for the algebra $A$, it remains to
construct a map
$$
HH_\idot(A) \to HH_\idot(Q_\idot(A)).
$$
To do this, one applies obstruction theory and shows that the map
$\alpha:Q_\idot(A) \to A$ admits a splitting in the category
$\Dgalg(k)$. The homology of the DG algebra $Q_\idot(A)$ is given by
\begin{equation}\label{Q.st}
\HH_i(Q_\idot(A)) = A \otimes \St_i(k),
\end{equation}
where $\St_\idot(k)$ is the dual $k$-Steenrod algebra --- that is,
the dual to the algebra of $k$-linear cohomological operations in
cohomology with coefficients in $k$. We have $\St_0(k) \cong
\St_1(k) \cong k$, and $\St_i(k) = 0$ for $1 < k \leq 2p$. The map
$a:Q_\idot(A) \to A$ is an isomorphism in degree $0$. The splitting
is constructed degree-by-degree. In degree $1$, the obstruction to
splitting is exactly the same as the obstruction to lifting the
algebra $A/k$ to the ring $W_2(k)$. In any higher degree $i > 1$,
the obstruction lies in the Hochschild cohomology group $HH^{2+i}(A
\otimes \St_i(k))$, and this vanishes in the relevant range of
degrees by the assumption $HH^i(A) = 0$, $i \geq 2p-1$.

In the DG algebra case, the construction breaks down since
Lemma~\ref{tt.le} does not have a DG version. Thus one first has to
replace a DG algebra $A^\hdot$ with a cosimplicial algebra $\A$ by
the Dold-Kan equivalence, and then apply the above construction to
$\A$ ``pointwise''. It is at this point that one has to require
$A^i=0$ for $i < 0$.

\medskip

Although \cite{K} only provides a Cartier map for DG algebras
defined over a finite field $k$, the same technology should apply to
DG algebras over $W=W(k)$ with very little changes, so that for any
smooth DG algebra $A^\hdot/W(k)$ with $HH^i(A^\hdot)=0$ for $i \geq
2p-1$, one should be able to construct a canonical isomorphism
$$
\wt{\phi}:\wt{HP_\idot(A^\hdot)} \cong HP_\idot(A^\hdot).
$$
Equivalently, $HP_\idot(A^\hdot)$ should carry a filtered
Dieudonn\'e module structure (in other words, underlie a canonical
object of the periodic derived category $\DFDM^{per}(W)$). One also
should be able to check that if $A^\hdot$ is Morita-equivalent to a
smooth variety $X/W$, the comparison isomorphism \eqref{hp.dr} is
compatible with the filtered Dieudonn\'e module structures on both
sides. However, at present, none of this has been done.

We note that the problem with the regulator map in the $p$-adic
setting mentioned in the end of Section~\ref{fdm.1.sec} survives in
the non-commutative situation. Namely, the standard technology for
constructing the regulator map \eqref{non-comm.reg} (\cite[Section
8.4]{Lo}) involves considering the group algebras $k[G]$ for $G =
GL_n(A)$, for all $n \geq 1$. As $n$ goes to infinity, the
homological dimension of these group alegbras becomes arbitrarily
large, and the conditions of Theorem~\ref{car} cannot be satified.

\section{Generalities on stable homotopy.}

The appearance of the Steenrod algebra in \eqref{Q.st} suggests that
the whole story should be related to algebraic topology. This is
indeed so. To explain the relation, we need to recall some standard
facts on stable homotopy theory.

\subsection{Stable homotopy category and homology.}
Roughly speaking, the {\em stable homotopy category} $\Sthom$ is
obtained by inverting the suspension functor $\Sigma$ in the
category $\hhom$ of pointed CW complexes and homotopy classes of
maps between them. Objects of $\Sthom$ are called {\em spectra}. A
spectrum consists of a collection of pointed CW complexes $X_i$, $i
\geq 0$, and maps $\Sigma X_i \to X_{i+1}$ for all $i$ (in some
treatments, these data are required to satisfy additional technical
conditions). For the definitions of maps between spectra and
homotopies between such maps, we refer the reader to a number of
standard references, for example \cite{adams}. Any CW complex $X \in
\hhom$ defines its {\em suspension spectrum} $\Sigma^\infty X \in
\Sthom$ consisting of the suspensions $\Sigma^i X$. For any two CW
complexes $X$, $Y$, we have
$$
\Hom_{\Sthom}(\Sigma^\infty X,\Sigma^\infty Y) =
\lim_{\overset{i}{\to}}[\Sigma^iX,\Sigma^iY],
$$
where $[-,-]$ denotes the set of homotopy classes of maps.

Any complex of abelian groups $M_\idot$ defines a spectrum
$\EM(M_\idot)$ called the {\em Eilenberg-Maclane spectrum of
$M_\idot$}. This is functorial in $M_\idot$, so that for any
commutative ring $R$, we have a functor
$$
\EM:\D(R) \to \D(\Ab) \to \Sthom,
$$
where $\D(R)$ is the derived category of the category of
$R$-modules.  This functor has a left-adjoint $H(R):\Sthom \to
\D(R)$, known as {\em homology with coefficients in $R$}.

The category $\Sthom$ is a tensor triangulated category. Both
functors $\EM$ and $H(R)$ are triangulated. Moreover, the homology
functor $H(R)$ is a tensor functor -- for any two spectra $X,Y \in
\Sthom$ with smash-product $X \wedge Y$, there exists a functorial
isomorphism
$$
H(R)(X) \lotimes_R H(R)(Y) \cong H(R)(X \wedge Y).
$$
The adjoint Eilenberg-Maclane functor $\EM$ is pseudotensor -- we
have a natural map
$$
\EM(V_\idot) \wedge \EM(W_\idot) \to \EM(V_\idot \lotimes_R W_\idot)
$$
for any two objects $V_\idot,W_\idot \in \D(R)$. Thus for any
associative ring object $\A$ in $\Sthom$, its homology $H(R)(\A)$ is
a ring object in $\D(R)$, and conversely, for any associative ring
object $A_\idot \in \D(R)$, the Eilenberg-Maclane spectrum
$\EM(A_\idot)$ is a ring object in $\Sthom$.

In the homological setting, we know that the structure of a ``ring
object in $\D(R)$'' is too weak, and the right objects to consider
are DG algebras over $R$. To define an analogous notion for spectra
is non-trivial, since the traditional topological interpretation of
spectra does not behave too well as far as the products are
concerned. Fortunately, new models for $\Sthom$ have appeared more
recently, such as for example {\em $S$-modules} of \cite{EMM}, {\em
  orthogonal spectra} of \cite{MM}, or {\em symmetric spectra} of
\cite{HSS}. All these approaches give equivalent results; to be
precise, let us choose for example the last one. As shown in
\cite{HSS}, symmetric spectra form a symmetric monoidal category;
denote it by $\Symm$.  Then in this paper, a {\rm ring spectrum}
will denote a monoidal object in $\Symm$, and $\Stalg$ will denote
the category of ring spectra considered up to a homotopy equivalence
(formally, this is defined by putting a closed model structure on
the category of ring monoidal objects in $\Symm$ whose weak
equivalences are homotopy equivalences of the underlying symmetric
spectra). The homology functor $H(R)$ and the Eilenberg-Maclane
functor $\EM$ extend to functors
$$
H(R):\Stalg \to \Dgalg(R),\qquad \EM:\Dgalg(R) \to \Stalg.
$$
where as in Section~\ref{hc.sec}, $\Dgalg(R)$ is the category of DG
algebras over $R$ considered up to a quasiisomorphism.

\subsection{Equivariant categories.}\label{equiv.subs}

For any compact group $G$, a pointed ``$G$-CW complex'' is a pointed
CW complex $X$ equipped with a continuos action of $G$ such that the
fixed-point subset $X^g \subset X$ is a pointed subcomplex for any
$g \in G$. We will denote by $\hhom(G)$ the category of pointed
$G$-CW complexes and $G$-equivariant homotopy classes of
$G$-equivariant maps between them. We note that for any closed
subgroup $H \subset G$, sending $X$ to the fixed-point subspace $X^H
\subset X$ gives a well-defined functor
$$
\hhom(G) \to \hhom.
$$
This functor is representable in the following sense: for any $X \in
\hhom(G)$, we have a homotopy equivalence
\begin{equation}\label{fxp.sp}
X^H \cong \Maps_G([G/H]_+,X),
\end{equation}
where $\Maps_G(-,-)$ means the space of $G$-equivariant maps with
its natural topology, and $[G/H]_+$ is the pointed $G$-CW complex
obtained by adding a (disjoint) marked point to the quotient $G/H$
with the induced topology and $G$-action.

To define a stable version of the category $\hhom(G)$, one could
again simply invert the suspension functor. However, there is a more
interesting alternative: by definition, $n$-fold suspension
$\Sigma^n$ is the smash-product with an $n$-sphere, and in the
equivariant setting, one can allow the sphere to carry a non-trivial
$G$-action. The corresponding equivariant stable category has been
constructed in \cite{may1}; it is known as the {\em genuine
$G$-equivariant stable homotopy category} $\Sthom(G)$. To define it,
one needs to fix a real representation $U$ of the group $G$ which is
equipped with a $G$-invariant inner product and contains every
finite-dimensional inner-product representation countably many
times; this is called a ``complete $G$-universe''. Then a genuine
$G$-equivariant spectrum is a collection of $G$-CW complexes $X(V)$,
one for each finite-dimensional $G$-invariant inner-product subspace
$V \subset U$, and maps $S^W \wedge X(V) \to X(V \oplus W)$, one for
each inner-product $G$-invariant subspace $V \oplus W \subset U$,
where $S^V$ is the one-point compactification of the underlying
topological space of the representation $V$, with its natural
$G$-action. As in the non-equivariant case, $\Sthom(G)$ is a tensor
triangulated category. We have a natural suspension spectrum functor
$\Sigma^\infty:\hhom(G) \to \Sthom(G)$, and for any two objects $X,Y
\in \hhom(G)$, we have
$$
\Hom_{\Sthom(G)}(\Sigma^\infty X,\Sigma^\infty Y) = \lim_{\overset{V
    \subset U}{\to}}[S^V \wedge X, S^V \wedge Y]_G,
$$
where $[-,-]_G$ is the set of $G$-homotopy classes of
$G$-equivariant maps, and the limit is over all the
finite-dimensional $G$-invariant inner-product subspaces $V \subset
U$. The category $\Sthom(G)$ does depend on $U$, but this is not too
drastic: all complete $G$-universes are isomorphic, and for any
isomorphism $U \cong U'$ between complete $G$-universes, there is a
``change of universe'' functor which is an equivalence between the
corresponding versions of $\Sthom(G)$.

Forgetting the $G$-action gives a natural forgetful functor
$\Sthom(G) \to \Sthom$, and equipping a spectrum with a trivial
$G$-action gives an embedding $\Sthom \to \Sthom(G)$. Thus for any
$X \in \Sthom$ and $Y \in \Sthom(G)$, we have a functorial smash
product $X \wedge Y \in \Sthom(G)$. This has an adjoint: for any
$X,Y \in \Sthom(G)$, we have a natural spectrum $\Maps_G(X,Y) \in
\Sthom$ such that for any $Z \in \Sthom$, there is a functorial
isomorphism
$$
\Hom_{\Sthom(G)}(Z \wedge X,Y) \cong \Hom_{\Sthom}(Z,\Maps_G(X,Y)).
$$
For any closed subgroup $H \subset G$ and any $X \in \Sthom(G)$, one
can extend \eqref{fxp.sp} and define the fixed point spectrum $X^H$ by
the same formula,
\begin{equation}\label{psi}
X^H = \Maps_G(\Sigma^\infty[G/H]_+,X).
\end{equation}
However, this does not commute with the suspension spectrum functor
$\Sigma^\infty$. In \cite{may1}, a second fixed-points functor is
introduced, called the {\em geometric fixed points functor} and
denoted $\Phi^H$. It does commute with $\Sigma^\infty$, and also
commutes with smash products, so that there are functorial
isomorphisms
$$
\Phi^H(\Sigma^\infty X) \cong \Sigma^\infty X^H, \qquad \Phi^H(X
\wedge Y) \cong \Phi^H(X) \wedge \Phi^H(Y)
$$
for any $X,Y \in \Sthom(G)$. For any $X \in \Sthom(G)$, there exists
a canonical map
\begin{equation}\label{can.phi}
\can:X^H \to \Phi^H(X),
\end{equation}
functorial in $X$. Moreover, let $N_H \subset G$ be the normalizer
of the subgroup $H \subset G$, and let $W_H = N_H/H$ be the
quotient. Then $\Phi^H$ can be extended to a functor
$$
\wh{\Phi}^H:\Sthom(G) \to \Sthom(W_H),
$$
and the same is true for the usual fixed-points functor $X \mapsto
X^H$ of \eqref{psi}. The map $\can$ of \eqref{can.phi} then lifts to
a map of $W_H$-equivariant spectra.  Here if $\Sthom(G)$ is defined
on a complete $G$-universe $U$, then $\Sthom(W_H)$ should be defined
on the complete $W_H$-universe $U^H$. The functor $\wh{\Phi}^H$ has
a right-adjoint which is a fully faithful embedding $\Sthom(W_H) \to
\Sthom(G)$ (for example, if $H=G$, then this is the trivial
embedding $\Sthom \to \Sthom(G)$).

\subsection{Mackey functors.}\label{mackey.subs}

Assume from now on that the compact group $G$ is a finite group with
discrete topology. It is not difficult to extend the homology
functor $H(R)$ to a functor
$$
H(R):\Sthom(G) \to \D(G,R)
$$
with values in the derived category of $R[G]$-modules. However, this
version of equivariant homology looses a lot of information such as
fixed points. A more natural target for equivariant homology is the
category of the so-called {\em Mackey functors}. To define them, one
considers an additive category $\B_G$ whose objects are $G$-orbits
$G/H$ for all subgroups $H \subset G$, and whose $\Hom$-groups are
given by
\begin{equation}\label{b.g}
\begin{aligned}
\B^G([G/H_1],[G/H_2]) &=
\Hom_{\Sthom(G)}(\Sigma^\infty[G/H_1]_+,\Sigma^\infty[G/H_2]_+)=\\
&= \pi_0(\Maps_G(\Sigma^\infty[G/H_1]_+,\Sigma^\infty[G/H_2]_+)).
\end{aligned}
\end{equation}
An {\em $R$-valued $G$-Mackey functor} (\cite{dress}, \cite{lind},
\cite{tD}, \cite{may2}) is an additive functor from $\B_G$ to the
category of $R$-modules. The category of such functors is an abelian
category, denoted $\M(G,R)$. 

More explicitly, for any subgroups $H_1,H_2 \subset G$, one can
consider the groupoid $\Qq([G/H_1],[G/H_2])$ of diagrams $[G/H_1]
\gets S \to [G/H_2]$ of finite sets equipped with a $G$-action, and
isomorphisms between such diagrams. Then disjoint union turns these
groupoids into symmetric monoidal categories, the Cartesian product
turns the collection $\Qq(-,-)$ into a $2$-category with objects
$[G/H]$, and it seems very likely that the mapping spectra
$\Maps_G(\Sigma^\infty[G/H_1]_+,\Sigma^\infty[G/H_2]_+)$ are in fact
obtained from the classifying spaces $|\Qq([G/H_1],[G/H_2])|$ of
symmetric monoidal groupoids $\Qq([G/H_1],[G/H_2])$ by group
completion. At present, this has not been proved (\cite{M0});
however, the corresponding isomorphism is well-known at the level of
$\pi_0$: we have
$$
\pi_0(\Maps_G(\Sigma^\infty[G/H_1]_+,\Sigma^\infty[G/H_2]_+)) \cong
\pi_0(\Omega B |\Qq([G/H_1],[G/H_2])|),
$$
so that the groups $\B_G(-,-)$ are given by
\begin{equation}\label{b.g.bis}
  \B^G([G/H_1],[G/H_2])
  = \Z[\Iso(\Qq([G/H_1],[G/H_2]))]/\{[S_1 \coprod S_2] - [S_1] - [S_2]\},
\end{equation}
where $\Iso$ means the set of isomorphism classes of objects. 

\medskip

For any $X \in \Sthom(G)$, individual homology groups $H_i(R)(X)$
can be equipped with a natural structure of a Mackey functor in such
a way that $H_i(R)(X)([G/H]) \cong H_i(R)(X^H)$, $H \subset G$ (for
more details, see \cite{may2}). To collect these into a single
homology functor $H(R)$, one has to work out a natural derived
version of the abelian category $\M(G,R)$. This has been done
recently in \cite{Ka-ma}. Roughly speaking, instead of $\pi_0$ in
\eqref{b.g}, one should the chain homology complexes $C_\idot(-,\Z)$
of the corresponding spectra, and one should set
$$
\B^G_\idot([G/H_1],[G/H_2]) =
C_\idot(\Maps_G(\Sigma^\infty[G/H_1]_+,\Sigma^\infty[G/H_2]_+),\Z).
$$
In practice, one replaces this with complexes which compute the
homology of the spectra obtained by group completion from the
symmetric monoidal groupoids $\Qq([G/H_1],[G/H_2])$. This can be
computed explicitly, so that the complexes $\B^G_\idot(-,-)$
introduced in \cite[Section 3]{Ka-ma} are given by an explicit
formula, and spectra are not mentioned at all. One then shows that
the collection $\B^G_\idot(-,-)$ is an $A_\infty$-category in a
natural way, and one defines the triangulated category $\DM(G,R)$ of
{\em derived $R$-valued $G$-Mackey functors} as the derived category
of $A_\infty$-functors from $\B^G_\idot$ to the category of
complexes of $R$-modules.

In general, the category $\DM(G,R)$ turns out to be different from
the derived category $\D(\M(G,R))$ (although both contain the
abelian category $\M(G,R)$ as a full subcategory). On the level of
slogans, one can hope that the category $\DM(G,R)$ is the ``brave
new product'' of the category $\Sthom(G)$ and the derived category
$\D(R)$ of $R$-modules, taken over the non-equivariant stable
homotopy category $\Sthom$, so that we have a diagram
$$
\begin{CD}
\D(\M(G,R)) @>>> \DM(G,R) @>>> \Sthom(G)\\
@. @VVV @VVV\\
@. \D(R) @>>> \Sthom,
\end{CD}
$$
where the square is Cartesian in some ``brave new'' sense. On a more
mundane level, it is expected that the triangulated category
$\DM(G,R)$ reflects the structure of the category $\Sthom(G)$ in the
following way.
\begin{enumerate}
\item There exists a symmetric tensor product $- \otimes -$ on the
  triangulated category $\DM(G,R)$, and for any subgroup $H \subset
  G$, we have natural triangulated fixed-point functors
  $\Phi^H,\Psi^H:\DM(G,R) \to \D(R)$.
\item There exists a natural triangulated equivariant homology
  functor
$$
H_G(R):\Sthom(G) \to \DM(G,R)
$$
and natural functorial isomorphisms
$$
\begin{aligned}
\Phi^H(H_G(R)(X)) &\cong H(R)(\Phi^H(X)),\\
\Psi^H(H_G(R)(X)) &\cong H(R)(X^H),\\
H_G(X \wedge Y) &\cong H_G(X) \otimes H_G(Y)
\end{aligned}
$$
for any $X,Y \in \Sthom(G)$, $H \subset G$.
\end{enumerate}
In fact, most of these statements has been proved in \cite{Ka-ma},
although only for the so-called ``Spanier-Whitehead category'', the
full triangulated subcategory in $\Sthom(G)$ spanned by the
suspension spectra of finite CW complexes (the only thing not proved
is the compatibility $\Psi^H(H_G(R)X) \cong H(R)(X^H)$ which
requires one to leave the Spanier-Whitehead category). It has been
also shown in \cite{Ka-ma} that as in the case of spectra, the fixed
point functor $\Phi^H$ extends to a functor
\begin{equation}\label{phi.h}
\wh{\Phi}^H:\DM(G,R) \to \DM(W_H,R)
\end{equation}
with a fully faithful right-adjoint. These fixed-points functors
allow one to give a very explicit description of the category
$\DM(G,R)$. Namely, let $I(G)$ be the set of conjugacy classes of
subgroups in $G$, and for any $c \in I(G)$, let
$$
\DM_c(G,R) \subset \DM(G,R)
$$
be the full subcategory of such $M \in \DM(G,R)$ that $\Phi^H(M) =
0$ unless $H \subset G$ is in the class $c$.

\begin{prop}[\cite{Ka-ma}]
For any $c \in I(G)$, $\DM_c(G,R) \subset \DM(G,R)$ is an admissible
triangulated subcategory, and for any subgroup $H \subset G$ is a
subgroup in the class $c$, the functor $\wh{\Phi}^H$ of \eqref{phi.h}
induces an equivalence
$$
\wh{\phi}^H:\DM_c(G,R) \cong \D(W_H,R).
$$
\end{prop}

Moreover, equip $I(G)$ with the partial order given by
inclusion. Then it has been shown in \cite{Ka-ma} that unless $c
\leq c'$, $\DM_c(G,R)$ is left-orthogonal to $\DM_{c'}(G,R)$, so
that $\DM_c(G,R)$, $c \in I(G)$ form a semiorhtogonal decomposition
of the triangulated category $\DM(G,R)$ indexed by the partially
ordered set $I(G)$ (for generalities on semiorthogonal
decompositions, see \cite{boka}). To describe the gluing data
between the pieces of this semiorthogonal decomposition, one
introduces the following.

\begin{defn}
  Assume given a finite group $G$ and a module $V$ over $R[G]$. The
  {\em maximal Tate cohomology} $\vH^\hdot_{max}(G,V)$ is given by
$$
\vH^\hdot_{max}(G,V) = \RHom^\hdot_{\D^b(G/R)/\Ind}(R,V),
$$
where $\RHom^\hdot$ is computed in the quotient $\D^b(G,R)/\Ind$ of
the bounded derived category $\D^b(G,R)$ by the full saturated
triangulated subcategory $\Ind \subset \D^b(G,R)$ spanned by
representrations $\Ind_G^H(W)$ induced from a representation $W$ of
a subgroup $H \subset G$, $H \neq G$.
\end{defn}

Then for any two subgroups $H \subset H' \subset G$ with conjucacy
classes $c,c' \in I$, $c \leq c'$, the gluing functor between
$\DM_c(G,R)$ and $\DM_{c'}(G,R)$ is expressed in terms of maximal
Tate cohomology of the group $W_H$ and its various subgroups.

This description turns out to be very effective because maximal Tate
cohomology often vanishes. For example, if the order of the group
$G$ is invertible in $R$, $\vH^\hdot_{max}(G,V)=0$ for any
$R[G]$-module $V$, and the category $\DM(G,R)$ becomes simply the
direct sum of the categories $\DM_c(G,R) \cong \D(W_H,R)$ (for the
abelian category $\M(G,R)$, a similar decomposition theorem has been
proved some time ago by J. Thevenaz \cite{Th1}). On the other hand,
if $R$ is arbitrary but the group $G = \Z/n\Z$ is cyclic, then
$\vH^\hdot_{max}(G,V)=0$ for any $V$ unless $n = p$ is prime, in
which case $\vH^\hdot_{max}(G,V)$ reduces to the usual Tate
cohomology $\vH^\hdot(G,V)$.

\section{Cyclotomic traces.}

Returning to the setting of Theorem~\ref{car}, we can now explain
the appearance of the Steenrod algebra in \eqref{Q.st}: up to a
quasiisomorphism, the DG algebra $Q_\idot(A)$ of \eqref{Q.st} is in
fact given by
$$
Q_\idot(A) = H(k)(\EM(A))^{k^*},
$$
where the $k^*$-invariants are taken with respect to the natural
action of the multiplicative group $k^*$ of the finite field $k$
induced by its action on $k$.

In particular, this shows that it is not necessary to use dimension
arguments to construct a splitting $A \to Q_\idot(A)$ of the
augmentation map $Q_\idot(A) \to A$. For example, if we are given a
ring spectrum $\A$ with homology DG algebra $A^\hdot = H(k)(\A)$,
then a canonical map
\begin{equation}\label{sp.spl}
A^\hdot = H(k)(\A) \to H(k)(\EM(H(k)(\A)))
\end{equation}
exists simply by adjunction, and being canonical, it is in
particular $k^*$-invariant. Thus for any DG algebra of the form
$A^\hdot = H(k)(\A)$, the same procedure as in the proof of
Theorem~\ref{car} allows one to construct a Cartier map. However, in
this case one can do much more -- namely, one can compare the
homological story with the theory of {\em cyclotomic traces} and
{\em topological cyclic homology} known in algebraic topology. Let
us briefly recall the setup (we mostly follow the very clear and
concise exposition in \cite{HM}).

\subsection{Topological cyclic homology.}

For any unital associative algebra $A$ over a ring $k$, the
Hochschild homology complex $CH_\idot(A)$ of Section~\ref{hc.sec} is
in fact the standard complex of a simplicial $k$-module $A_\# \in
\Delta^{opp}k\amod$. {\em Topological Hochschild homology} is a
version of this construction for ring spectra. It was originally
introduced by B\"okstedt \cite{bok} long before the invention of
symmetric spectra, and used the technology of ``functors with a
smash product''. In the language of symmetric spectra, one starts
with a unital associative ring spectrum $\A$, and one defines a
simplicial spectrum $\A_\#$ by exactly the same formula as in the
algebra case. The terms of $\A_\#$ are the iterated smash products
$\A \wedge \dots \wedge \A$, and the face and degeneracy maps are
obtained from the multiplication and the unit map in $\A$. Then one
sets
$$
\THH(\A) = \hocolim_{\Delta^{opp}}\A_\#.
$$
As in the algebra case, this spectrum is equipped with a canonical
$S^1$-action, but in the topological setting this means much more:
one shows that $\THH(\A)$ actually underlies a canonical
$S^1$-equivariant spectum $\THH(\A) \in \Sthom(S^1)$.

However, this is not the end of the story. Note that the finite
subgroups in $S^1$ are the cyclic groups $C_n = \Z/n\Z \subset S^1$
numbered by integers $n \geq 1$, and for every $n$, we have $S^1/C_n
\cong S^1$. Fix a system of such isomorphisms which are compatible
with the embeddings $C_n \subset C_{nm} \subset S^1$, $n,m \geq 1$,
and fix a compatible system of isomorphisms $U^{C_n} \cong U$, where
$U$ is the complete $S^1$-universe used to define
$\Sthom(S^1)$. Then the following notion has been introduced in
\cite{BM}.

\begin{defn}\label{cyclo}
  \mbox{} A {\em cyclotomic structure} on an $S^1$-equivariant
  spectrum $T$ is given by a collection of $S^1$-equivariant
  homotopy equivalences
$$
r_n:\wh{\Phi}^{C_n}T \cong T,
$$
one for each finite subgroup $C_n \subset S^1$, such that $r_1= \id$
and $r_n \circ r_m = r_{nm}$ for any two integer $n,m > 1$.
\end{defn}

\begin{remark}\label{cyclo.rem}
  Here it is tacitly assumed that one works with specific model of
  equivariant spectra, so that a spectrum means more than just an
  object of the triangulated category $\Sthom(S^1)$; moreover, the
  functors $\wh{\Phi}^{C_n}$ are composed with the change of
  universe functors so that we can treat them as endofunctors of
  $\Sthom(S^1)$. Please refer to \cite{BM} or \cite{HM} for exact
  definitions.
\end{remark}

\begin{exa}\label{loop.exa}
  Assume given a CW complex $X$, and let $LX = \Maps(S^1,X)$ be its
  free loop space. Then for any finite subgroup $C \subset S^1$, the
  isomorphism $S^1 \cong S^1/C$ induces a homeomorphism
$$
\Maps(S^1,X)^C = \Maps(S^1/C,X) \cong \Maps(S^1,X),
$$
and these homeomorphism provide a canonical cyclotomic structure on
the suspension spectrum $\Sigma^\infty LX$.
\end{exa}

For any $S^1$-equivariant spectrum $T$ and a pair of integers $r,s >
1$, one has a natural non-equivariant map
$$
F_{r,s}:T^{C_{rs}} \to T^{C_r}.
$$
On the other hand, assume that $T$ is equipped with a cyclotomic
structure. Then we have a natural map
$$
\begin{CD}
R_{r,s}:T^{C_{rs}} \cong (T^{C_s})^{C_r} @>{\can}>>
(\wh{\Phi}^{C_s}T)^{C_r} @>{r_s}>> T^{C_r},
\end{CD}
$$
where $\can$ is the canonical map \eqref{can.phi}, and $r_s$ comes
from the cyclotomic structure on $T$.  To pack together the maps
$F_{r,s}$, $R_{r,s}$, it is convenient to introduce a small category
$\I$ whose objects are all integers $n \geq 1$, and whose maps are
generated by two maps $F_r,R_r:n \to m$ for each pair $m$, $n=rm$,
$r > 1$, subject to the relations $F_r \circ F_s = F_{rs}$, $R_r
\circ R_s = R_{rs}$, $F_r \circ R_s = R_s \circ F_r$. Then the maps
$T_{r,s}$, $F_{r,s}$ turn the collection $T^{C_n}$, $n \geq 1$ into
a functor $\wt{T}$ from $\I$ to the category of spectra.

\begin{defn}
The {\em topological cyclic homology} $\TC(T)$ of a cyclotomic
spectrum $T$ is given by
$$
\TC(T) = \holim_{\I} \wt{T}.
$$
\end{defn}

Given a ring spectrum $\A$, B\"okstedt and Madsen equip the
$S^1$-equivariant spectrum $THH_\idot(\A)$ with a canonical
cyclomomic structure. {\em Topological cyclic homology} $\TC(\A)$ is
then given by
$$
\TC(\A) = \TC(\THH(\A)).
$$
Further, they construct a canonical {\em cyclotomic trace map}
\begin{equation}\label{cyc.tr}
K(\A) \to \TC(\A)
\end{equation}
from the $K$-theory spectrum $K(\A)$ to the topological cyclic
homology spectrum.

The topological cyclic homology functor $\TC(\A)$ and the cyclotomic
trace were actually introduced by B\"okstedt, Hsiang and Madsen in
\cite{BHM}; the more convenient formulation using cyclotomic spectra
appeared slightly later in \cite{BM}. Starting with \cite{BHM}, it
has been proved in many cases that the cyclotomic trace map becomes
a homotopy equivalence after taking profinite completions of both
sides of \eqref{cyc.tr}. Moreover, in \cite{mac} MacCarthy
generalized Goodwillie's Theorem and proved that after pro-$p$
completion at any prime $p$, the cyclotomic trace gives an
equivalence of the relative groups $\wh{K(\A,I)}_p \cong
\wh{\TC(\A,I)}_p$, where $I \subset \A$ is a nilpotent ideal.

\subsection{Cyclotomic complexes.}

To define a homological analog of cyclotomic spectra, one needs to
replace $S^1$-equivariant spectra with derived Mackey functors. The
machinery of \cite{Ka-ma} does not apply directly to non-discrete
groups, since this would require treating the groupoids $\Qq(-,-)$
of Subsection~\ref{equiv.subs} as topological groupoids. However,
for finite subgroups $C_1,C_2 \subset S^1$, the category
$\Qq([S^1/C_1],[S^1/C_2])$ is still discrete. Thus one can define a
restricted version of derived $S^1$-Mackey functors by discarding
the only infinite closed subgroup in $S^1$ (which is $S^1$
itself). This is done in \cite{cyclo}. The category $\DML(R)$ of
{\em $R$-valued cyclic Mackey functors} introduced in that paper has
the following features.
\begin{enumerate}
\item For every proper finite subgroup $C = C_n \subset S^1$, $n >
  1$, there is a fixed-point functor $\wh{\Phi}_n:\DML(R) \to
  \DML(R)$ whose right-adjoint functor $\wh{\iota}_n:\DML(R) \to
  \DML(R)$ is a full embedding. Moreover, there are canonical
  isomorphisms $\wh{\Phi}_n \circ \wh{\Phi}_m \cong \wh{\Phi}_{mn}$.
\item Let $\D_{S^1}(R)$ be the equivariant derived category of
  Section~\ref{hdr.sec}. Then there is a full embedding
  $\iota_1:\D_{S^1}(R) \to \DML(R)$ with a left-adjoint
  $\Phi_1:\DML(R) \to \D_{S^1}(R)$.
\item The images $\DML_n(R)$ of the full embeddings $\iota_n =
  \wh{\iota}_N \circ \iota_1:\D_{S^1}(R) \to \DML(R)$, $n \geq 1$,
  generate the triangulated category $\DML(R)$, and $\DML_n(R)
  \subset \DML(R)$ is left-orthogonal to $\DML_m(R) \subset \DML(R)$
  unless $n=mr$ for some integer $r \geq 1$.
\end{enumerate}
Thus as in the finite group case of \cite{Ka-ma}, the subcategories
$\DML_n(R) \subset \DML(R)$ form a semiorthogonal decomposition of
the category $\DML(R)$. The gluing data between $\DML_{mr}(R)$ and
$\DML_r(R)$ can be expressed in terms of the maximal Tate cohomology
$\vH^\hdot_{max}(C_m,-)$ of the cyclic group $C_m = \Z/m\Z$. For any
$n \geq 1$, let $\overline{\Phi}_n:\DML(R) \to \D(R)$ be the
composition of the left-adjoint $\Phi_n = \Phi_1 \circ \wh{\Phi}_n$
to $\iota_n$ and the forgetful functor $\D_{S^1}(R) \to \D(R)$; then
the functors $\overline{\Phi}_n$ play the role of fixed points
functors $\Phi^H$. There are also functors $\Psi_n:\DML(R) \to
\D_{S^1}(R)$ analogous to the functors $\Psi^H$. The homology
functor $H(R)$ extends to a functor
$$
H_{S^1}(R):\Sthom(S^1) \to \DML(R),
$$
and we have functorial isomorphisms
$$
\overline{\Phi}_n(H_{S^1}(R)(T)) \cong H(R)(\Phi^{C_n}(T)), \qquad
\Psi_n(H_{S^1}(R)(T)) \cong H(R)(T^{C_n})
$$
for every $n \geq 1$ and every $T \in \Sthom(S^1)$.

Another category defined in \cite{cyclo} is a triangulated category
$\DLR(R)$ of {\em $R$-valued cyclotomic complexes}. Essentially, a
cyclotomic complex $M_\idot \in \DLR(R)$ is a cyclic Mackey functor
$M_\idot$ equipped with a system of compatible quasiisomorphisms
$$
\wh{\Phi}_nM_\idot \cong M_\idot,
$$
as in Definition~\ref{cyclo} (although as in Remark~\ref{cyclo.rem},
the precise definition is different for technical reasons). The
homology functor $H_{S^1}(R):\Sthom(S^1) \to \DML(R)$ extends to a
functor from the category of cyclotomic spectra to the category
$\DLR(R)$. Moreover, all the constructions used in the definition of
topological cyclic homology make sense for cyclotomic complexes, so
that one has a natural functor
$$
\TC:\DLR(R) \to \D(R)
$$
and a functorial isomorphism
\begin{equation}\label{tc.tc}
\TC(H_{S^1}(R)(T)) \cong H(R)(\TC(T))
\end{equation}
for every cyclotomic spectrum $T$.

\subsection{Comparison theorem.}

We can now formulate the comparison theorem relating Dieudonn\'e
modules and cyclotomic complexes. We introduce the following
definition.

\begin{defn}\label{gfdm}
A {\em generalized filtered Dieudonn\'e module} $M$ over a
commutative ring $R$ is an $R$-module $M$ equipped with a decreasing
filtration $F^\hdot M$ and a collection of maps
$$
\phi^p_{i,j}:F^iM \to M/p^j,
$$
one for every integers $i$, $j$, $j \geq 1$, and a prime $p$, such
that
$$
\phi^p_{i,j+1} = \phi^p_{i,j} \mod p^j, \quad
\phi^p_{i,j}|_{F^{i+1}M} = p\phi^p_{i,j}.
$$
\end{defn}

For any integer $i$, we define the generalized filtered Dieudonn\'e
module $R(i)$ as $R$ with the filtration $F^iR(i)=R$,
$F^{i+1}R(i)=0$, and $\phi^p_{i,j} = p^i\id$ for any $p$ and $j$.
Generalized filtered Dieudonn\'e modules in the sense of
Defintion~\ref{gfdm} do not form an abelian category; however, by
inverting the filtered quasiisomorphisms, we can still construct the
derived category $\DFDM_{g}(R)$ and its twisted $2$-periodic version
$\DFDM_{g}^{per}(R)$. 

Definition~\ref{gfdm} generalizes \eqref{fdm.defn} in that it
collects together the data for all primes $p$. Note, however, that
one can rephrase Definition~\ref{gfdm} by putting together all the
maps $\phi^p_{i,j}$, $j \geq 1$, into a single map
$$
\wh{\phi}^p_i:F^iM \to \wh{(M)}_p
$$
into the pro-$p$ completion $\wh{(M)}_p$ of the module
$M$. Then if $R = \Z_p$ and $M$ is finitely generated over $\Z_p$,
we have
$$
\wh{(M)}_p \cong M, \quad \wh{(M)}_l = 0 \text{ for } l \neq p,
$$
so that for such an $M$, the extra data imposed onto $M$ in
Definition~\ref{gfdm} and in Definition~\ref{fdm.defn} are the same.
In general, for any prime $p$, we have a fully faithful embedding
$$
\wt{\DFDM}(\Z_p) \subset \DFDM_{g}(\Z),
$$
where $\wt{\DFDM}(\Z_p)$ is as in Section~\ref{fdm.1.sec}, and
similarly for the periodic categories. The essential images of these
embeddings are spanned by complexes which are pro-$p$ complete as
complexes of abelian groups. Note, however, that what appears here
are {\em weak} filtered Dieudonn\'e modules. The requirement that
the map \eqref{wt.phi} is a quasiisomorphism can be additionally
imposed at each individual prime $p$; I do not know whether it is
useful to impose it in the universal category $\DFDM_{g}(\Z)$.

Here is then the main comparison theorem of \cite{cyclo}.

\begin{theorem}[{{\cite[Section 5]{cyclo}}}]\label{cyclo.thm}
For any commutative ring $R$, there is a canonical equivalence of
categories
$$
\DFDM_{g}^{per}(R) \cong \DLR(R).
$$
\end{theorem}

Thus the category $\DLR(R)$ of cyclotomic complexes over $R$ admits
an extremely simple linear-algebraic description. Roughly speaking,
the reason for this is the vanishing of maximal Tate cohomology
$\vH^\hdot(\Z/n\Z,-)$ for non-prime $n$ mentioned at the end of
Subsection~\ref{mackey.subs}. Due to this vanishing, the only
non-trivial gluing between the pieces $\DML_m(R)$, $\DML_n(R)$ of
the semiorthogonal decomposition of the category $\DML(R)$ of cyclic
Mackey functors occurs when $n = mp$ for a prime $p$ (and this
gluing is described by the Tate cohomology of the group
$\Z/p\Z$). The gluing data provide the maps $\wh{\phi}^p_i$ in the
equivalence of Theorem~\ref{cyclo.thm}; the periodic filtered
complex comes from the equivalence $\D_{S^1}(R) \cong \DF^{per}(R)$
of Lemma~\ref{trivial.lemma}. These are the main ideas of the proof.

Moreover, there is a second comparison theorem which expresses
topological cyclic homology in terms of generalized Dieudonn\'e
modules.

\begin{theorem}[{{\cite[Section 6]{cyclo}}}]\label{tc.thm}
  Under the equivalence of Theorem~\ref{cyclo.thm}, there is
  functorial isomorphism
\begin{equation}\label{tc.eq}
\wh{\TC(M_\idot)}_f \cong \wh{\RHom^\hdot(R,M_\idot)}_f
\end{equation}
for any $M_\idot \in \DLR(R)$, where $R = R(0) \in
\DFDM^{per}_{g}(R)$ is the trivial generalized filtered Dieudonn\'e
module, and $\wh{(-)}_f$ stands for profinite completion.
\end{theorem}

\begin{remark}
  Both $\TC(-)$ and $\RHom^\hdot(R,-)$ commute with profinite
  completions, so that if $M_\idot$ itself is profinitely complete,
  the completions in \eqref{tc.eq} can be dropped. In general, it is
  better to keep the completion; to obtain an isomorphism in the
  general case, one should, roughly speaking, replace $T$ with the
  homotopy fixed points $T^{hS^1}$ in the definition of topological
  cyclic homology $\TC(T)$.
\end{remark}

\begin{remark}
  It is not unreasonable to hope that Theorem~\ref{tc.thm} has a
  topological analog: one can define a triangulated category of
  cyclotomic spectra which is enriched over $\Sthom$, and then for
  any profinitely complete cyclotomic spectrum $T$, we have a
  natural homotopy equivalence
$$
\TC(T) \cong \Maps(\Ss,T),
$$
where $\Maps(-,-)$ is the mapping spectrum in the cyclotomic
category, and $\Ss = \Sigma^\infty\ppt$ is the sphere spectrum with
the trivial cyclotomic structure (obtained as in
Example~\ref{loop.exa}). This would give a conceptual replacement
of the somewhat {\em ad hoc} definition of the functor $\TC$.
\end{remark}

\subsection{Back to ring spectra.} 

Return now to our original situation: we have a ring spectrum $\A
\in \Stalg$, and the DG algebra $A_\idot = H(W)(\A)$ is obtained as
its homology with coefficients in the Witt vector ring $W= W(k)$ of
a finite field $k$. Assume for simplicity that $k = \Z/p\Z$ is a
prime field, so that $W=\Z_p$.

Then on one hand, we have the cyclotomic spectrum $\THH(\A)$ of
\cite{BM}, and since the homology functor $H(\Z_p)$ commutes with
tensor products, we have a quasiisomorphism
$$
H(\Z_p)(\THH(\A)) \cong CH_\idot(A_\idot).
$$
But the left-hand side underlies a cyclotomic complex, and by
Theorem~\ref{cyclo.thm}, this is equivalent to saying that it has a
structure of a generalized filtered Dieudonn\'e module. And on the
other hand, $CH_\idot(A_\idot)$ has a Dieudonn\'e module structure
induced by the splitting map \eqref{sp.spl}. We expect that the two
structures coincide (although at present, this has not been checked).

Moreover, the functor $\TC$ commutes with $H(\Z_p)$ by
\eqref{tc.tc}, and Theorem~\ref{tc.thm} shows that we have
$$
\begin{aligned}
H(\Z_p)(\TC(\A)) &\cong H(\Z_p)(\TC(\THH(\A))) \cong
\TC(CH_\idot(A_\idot))\\
&\cong \RHom^\hdot(\Z_p,CH_\idot(A)),
\end{aligned}
$$
where $\RHom^\hdot(-,-)$ is taken in the category
$\DFDM^{per}(\Z_p)$ of filtered Dieudonn\'e modules. In other words:
\begin{itemize}
\item The homology functor $H(\Z_p)$ sends topological cyclic
  homology into syntomic periodic cyclic homology.
\end{itemize}
This principle can be used to study further the regulator map for
syntomic homology. Namely, applying $H(\Z_p)$ to the cyclotomic
trace map \eqref{cyc.tr}, we obtain a functorial map
$$
H(\Z_p)(K_\idot(\A)) \to H(\Z_p)(\TC(\A)),
$$
and the right-hand side is the target of the desired regulator map
for the DG algebra $A_\idot$. The desired source of this map is
$K_\idot(A_\idot) \cong K_\idot(H(\Z_p)(\A))$. Thus the question of
existence of the syntomic regulator maps reduces to a problem in
algebraic $K$-theory: describe the relation between the homology of
the $K$-theory of a ring spectrum, and the $K$-theory of its
homology.

\medskip

To finish the Section, let us explain how things work in a very
simple particular case. Assume given a CW complex $X$, and let $\A =
\Sigma^\infty \Omega X$, the suspension spectrum of the based loop
space $\Omega X$. Then since $\Omega X$ is a topological monoid,
$\A$ is a ring spectrum. The DG algebra $A_\idot = H(\Z_p)(\A)$ is
given by $A_\idot = C_\idot(\Omega X,\Z_p)$, the singular chain
complex of the topological space $\Omega X$. It is known that in
this case, we have
$$
CH_\idot(A_\idot) \cong C_\idot(LX,\Z_p),
$$
the singular chain complex of the free loop space $LX$. Analogously,
we have $\THH(\A) \cong \Sigma^\infty LX$. The $S^1$-action on
$\THH(\A)$ and $CH_\idot(A_\idot)$ is induced by the loop rotation
action on $LX$. The cyclotomic structure on $\THH(\A)$ is that of
Example~\ref{loop.exa}. The corresponding Dieudonn\'e module
structure map $\phi$ on $CP_\idot(A_\idot)$ is induced by the
cyclotomic structure map $LX^{\Z/p\Z} \cong LX$ of the free loop
space $LX$. To compare this with the constructions of
Section~\ref{fdm.2.sec}, specialize even further and assume that
$\Omega X$ is discrete, so that $A_\idot$ is quasiisomorphic to an
algebra $A$ concentated in degree $0$. In this case $X \cong BG$ for
a discrete group $G$, and $A = \Z_p[G]$ is its group algebra. Then
the diagonal map $G \to G^p$ induces a map $A \to A^{\otimes p}$
which is a quasi-Frobenius map in the sense of
Section~\ref{fdm.2.sec}, thus induces another Dieudonn\'e module
structure on the filtered complex $CP_\idot(A)$. One checks easily
that the two structures coincide. For a general $X$, the Diedonn\'e
module structure on $CP_\idot(A_\idot)$ can also be described
explicitly in the same way as in Section~\ref{fdm.2.sec}, by using
the map
$$
A_\idot \to A_\idot^{\otimes p}
$$
induced by the diagonal map $\Omega X \to (\Omega X)^p$ in place of
the quasi-Frobenius map.

\section{Hodge structures.}

In the archimedian setting of \thetag{i} of Section~\ref{mm.sec},
much less is known about periodic cyclic homology than in the
non-archimedian setting of \thetag{ii}. One starts with a smooth
proper DG algebra $A^\hdot$ over $\C$ and considers its periodic
cyclic homology complex $CP_\idot(A^\hdot)$ with its Hodge
filtration. In order to equip $HP_\idot(A^\hdot)$ with an $\R$-Hodge
structure, one needs to define a weight filtration $W_\idot
CP_\idot(A^\hdot)$ and a complex conjugation isomorphism
$\overline{\phantom{m}}:CP_\idot(A^\hdot) \to
\overline{CP_\idot(A^\hdot)}$. The gradings in the isomorphism
\eqref{hp.dr} suggest that $W_\idot$ should be simply the canonical
filtration of the complex $CP_\idot(A^\hdot)$. However, the complex
conjugation is a complete mystery. There is only one approach known
at present, albeit a very indirect and highly conjectural one; the
goal of this section is to describe it. I have learned all this
material from B. To\"en and/or M. Kontsevich -- it is only the
mistakes here that are mine.

The so-called {\em $\D^-$-stacks} introduced by B. To\"en and
G. Vezzosi in \cite{ToVe2} generalize both Artin stacks and DG
schemes and form the subject of what is now known as ``derived
algebraic geometry''; a very nice overview is avaiable in
\cite{toen-overview}. Very approximately, a $\D^-$-stack over a ring
$k$ is a functor
$$
\M:\Delta^{opp}\Comm(k) \to \Delta^{opp}\Sets
$$
from the category of simplicial commutative algebras over $k$ to the
category of simplicial sets. This functor should satisfy some
descent-type conditions, and all such functors are considered up to
an appropriately defined homotopy equivalence (made sense of by the
technology of closed model structures). This generalizes the
Grothendieck approach to schemes which treats a scheme over $k$ as
its functor of points -- a sheaf of sets on the opposite
$\Comm(k)^{opp}$ to the category of commutative algebras over
$k$. The category $\Comm(k)$ is naturally embedded in
$\Delta^{opp}\Comm(k)$ as the subcategory of constant simplicial
objects, and resticting a $\D^-$-stack $\M$ to $\Comm(k) \subset
\Delta^{opp}\Comm(k)$ gives an $\infty$-stack in the sense of
Simpson \cite{Si} (this is called the {\em truncation} of $\M$).

If $k$ contains $\Q$, one may replace simplicial commutative
algebras with commutative DG algebras $R_\idot$ over $k$ placed in
non-negative homological degrees, $R_i = 0$ for $i < 0$. If we
denote the category of such DG algebras by $\Dgcomm^-(k)$, then a
$\D^-$-stack is a functor
$$
\M:\Dgcomm^-(k) \to \Delta^{opp}\Sets,
$$
again satisfying some conditions, and considered up to a homotopy
equivalence. The category of $\D^-$-stacks over $k$ is denoted
$\Dst(k)$. For every DG algebra $R_\idot \in \Dgcomm^-(k)$, its
{\em derived spectrum} $\Rspec(R_\idot) \in \Dst(k)$ sends a DG
algebra $R'_\idot \in \Dgcomm^-(k)$ to the simplicial set of maps
from $R_\idot$ to $R'_\idot$, with the simplicial structure induced
by the model structure on the category $\Dgcomm^-(k)$. We thus
obtain a Yoneda-type embedding
$$
\Rspec:\Dgcomm^-(k)^{opp} \to \Dst(k).
$$
For any DG algebra $R_\idot \in \Dgcomm^-(k)$, its de Rham
cohomology comple $\Omega^\hdot(R_\idot)$ is defined in the obvious
way; $\Omega^\hdot(-)$ gives a functor
$$
\Omega^\hdot:\Dgcomm^-(k) \to \Spaces_\Q
$$
from $\Dgcomm^-(k)$ to the category $\Spaces_\Q$ of rational homotopy
types in the sense of Quillen \cite{Q}. By the standard Kan
extension machinery, $\Omega^\hdot$ extentds to a de Rham
realization functor
$$
\Omega^\hdot:\Dst(k) \to \Spaces_\Q.
$$
Alternatively, one can take the $0$-th homology algebra
$H_0(R_\idot)$ and consider its cristalline cohomology; this gives a
DG algebra quasiisomorphic to $\Omega^\hdot(R_\idot)$ (the higher
homology groups behave as nilpotent extensions and do not contribute
to cohomology). This shows that the de Rham realization
$\Omega^\hdot(\M)$ of a $\D^-$-stack $\M \in \Dst(k)$ only depends
on its truncation.

Moreover, for $\D^-$-stacks satifying a certain finiteness condition
(``locally geometric'' and ``locally finitely presented'' in the
sense of \cite{ToVa}), instead of considering de Rham cohomology,
one can take the underlying topological spaces $\Top(\M(R))$ of the
simplicial complex algebraic varieties $\M(R)$, $R \in \Comm(k)$; by
Kan extension, this gives a topological realization functor
$$
\Top:\Dst(k) \to \Spaces
$$
into the category of topological spaces. By the standard comparison
theorems, $\Top(\M)$ and $\Omega^\hdot(\M)$ represented the same
rational homotopy type.

Now, it has been proved in \cite{ToVa} that for any associative
unital DG algebra $A^\hdot$ over $k$, there exists a $\D^-$-stack
$\M(A^\hdot)$ classifying ``finite-dimensional DG modules over
$A^\hdot$''. By definition, for any commutative DG algebra $R_\idot
\in \Dgcomm^-(k)$, the simplicial set $\M(A^\hdot)(R_\idot)$ is
given by
\begin{itemize}
\item $\M(A^\hdot)(R_\idot)$ is the nerve of the category
  $\Perf(A^\hdot,R_\idot)$ of DG modules over $A^\hdot \otimes
  R_\idot$ which are perfect over $R_\idot$, and quasiisomorphisms
  between such DG modules.
\end{itemize}
To\"en and Vaqui\'e prove that this indeed defines a $\D^-$-stack.
Moreover, they prove that if $A^\hdot$ satisfies certain finiteness
conditions, the $\D^-$-stack $\M(A^\hdot)$ is locally geometric and
locally finitely presented.

In particular, a smooth and proper DG algebra $A^\hdot \in
\Dgalg(k)$ satisfies the finiteness conditions needed for
\cite{ToVa}, so that there exists a locally geometric and locally
finitely presented $\D^-$-stack $\M(A^\hdot)$. Consider its de Rham
realization $\Omega^\hdot(\M(A^\hdot))$. For any $R_\idot \in
\Dgcomm^-(k)$, the category $\Perf(A^\hdot,R_\idot)$ is a symmetric
monoidal category with respect to the direct sum, so that the
realization $\Top(\M(A^\hdot))$ is automatically an
$E_\infty$-space.

\begin{lemma}[To\"en]\label{to.1}
  The $E_\infty$-space $\Top(\M(A^\hdot))$ is group-like.
\end{lemma}

\proof[Sketch of a possible proof.] One has to show that
$\pi_0(\Top(\M(A^\hdot)))$ is not only a commutative monoid but also
an abelian group. A point in $\Top(\M(A^\hdot))$ is represented by a
DG module $M_\idot$ over $A^\hdot$ which is perfect over $k$. One
observes that $M_\idot \oplus M_\idot[1]$ can be deformed to a
acyclic DG module; thus the sum of points represented by $M_\idot$
and $M_\idot[1]$ lies in connected component of $0$ in
$\Top(\M(A^\hdot))$.
\endproof

Thus for any smooth and proper DG algebra $A^\hdot \in \Dgalg(k)$,
the realization $\Top(\M(A^\hdot))$ is an infinite loop space, that
is, the $0$-th component of a spectrum.

\begin{defn}\label{K.st}
The {\em semi-topological $K$-theory} $K_\idot^{st}(A^\hdot)$ of a
smooth and proper DG algebra $A^\hdot$ is given by
$$
K_\idot^{st}(A^\hdot) = \pi_\idot(\Top(\M(A^\hdot))),
$$
the homotopy groups of the infinite loop space $\Top(\M(A^\hdot))$.
\end{defn}

If we are only interested in $K^{st}_\idot(A^\hdot) \otimes k$, we
may compute it using the de Rham model
$\Omega^\hdot(\M(A^\hdot))$. Then $K_\idot^{st}(\M(A^\hdot))$ is
exactly the complex of primitive elements with respect to the
natural cocommutative coalgebra structure on $\M(A^\hdot)$ induces
by the direct sum map
$$
\M(A^\hdot) \times \M(A^\hdot) \to \M(A^\hdot).
$$
Since $\Q \subset k$, and rationally, spectra are the same as
complexes of $\Q$-vector spaces, the groups $K^{st}_\idot(A^\hdot)
\otimes k$ are the only rational invariants one can extract from
the space $\M(A^\hdot)$.

Assume for the moment that $A^\hdot \in \Dgalg(k)$ is derived-Morita
equivalent to a smooth and proper algebraic variety $X/k$. Then one
can also consider the $\infty$-stack $\overline{\M}(X)$ of all
coherent sheaves on $X$; for any noetherian $R \in \Comm(k)$,
$\overline{\M}(X)(R)$ is by definition the nerve of the category of
coherent sheaves on $M \otimes R$ and isomorphisms between them. The
realization $\Top(\overline{\M}(X))$ is again an $E_\infty$-space,
no longer group-like. By definition, we have a natural map
$$
\overline{\M}(X) \to \M(A^\hdot),
$$
and the induced $E_\infty$-map of realizations. 

\begin{lemma}[To\"en]\label{to.2}
The natural $E_\infty$-map
\begin{equation}\label{bar.m}
\Top(\overline{\M}(X)) \to \Top(\M(A^\hdot))
\end{equation}
induces a homotopy equivalence between $\Top(\M(A^\hdot))$ and the
group completion of the $E_\infty$-space $\Top(\overline{\M}(X))$.
\end{lemma}

\proof[Sketch of a possible proof.] Since $\Top(\M(A^\hdot))$ is
group-like by Lemma~\ref{to.1}, it suffices to prove that the
delooping
$$
B\Top(\overline{\M}(X)) \to B\Top(\M(A^\hdot))
$$
of the $E_\infty$-map \eqref{bar.m} is a homotopy
equivalence. Delooping obviously commutes with geometric
realization, so that $B\Top(\M(A^\hdot))$ is the realization of the
$\D^-$-stack $B\M(A^\hdot)$, and similarly for
$B\Top(\overline{\M}(X))$. Instead of taking deloopings, we can
apply Waldhausen's $S$-construction. The resulting map
$$
S\overline{\M}(X) \to S\M(A^\hdot)
$$
is then an equivalence by Waldhausen's devissage theorem, so that it
suffices to prove that the natural map
$$
\Top(B\M(A^\hdot)) \to \Top(S\M(A^\hdot))
$$
is a homotopy equivalence, and similarly for $\overline{\M}(X)$. For
this, one argues as in Lemma~\ref{to.1}: since every filtered
complex can be canonically deformed to its associated graded
quotient, the terms $\Top(S_n\M(A^\hdot))$ of the $S$-construction
can be retracted to $n$-fold products $\Top(\M(A^\hdot) \times \dots
\times \M(A^\hdot))$, that is, the terms of the delooping
$\Top(B\M(A^\hdot))$, and similarly for $\overline{\M}(X)$.
\endproof

\begin{corr}\label{to.c}
The semitopological $K$-theory $K^{st}_\idot(\Q)$ is given by
$$
K^{st}_\idot(k) \cong \Z[\beta],
$$
the algebra of polynomials in one generator $\beta$ of degree $2$.
\end{corr}

\proof{} By Lemma~\ref{to.2}, computing $K^{st}_\idot(k)$ reduces to
studying the group completion of the realization
$$
\Top(\overline{\M}(\ppt)) \cong \coprod_n \Top([\ppt/GL_n]) \cong
\coprod_n BU_n,
$$
where $[\ppt/GL_n]$ is the Artin stack obtained as the quotient of
the point by the trivial action of the algebraic group $GL_n$.  This
group completion is well-known to be homotopy equivalent to the
classifying space $\Z \times BU$.
\endproof

\begin{remark} At present, Lemma~\ref{to.1} and Lemma~\ref{to.2} are
  unpublished, as well as Corollary~\ref{to.c}. The above sketches
  of proofs have been kindly explained to me by
  B. To\"en. Lemma~\ref{to.1} is slightly older, and it also appears
  for example in \cite{To.lec}.
\end{remark}

Now, since $k \supset \Q$ by our assumption, we have a well-defined
tensor product $M_\idot \otimes V_\idot$ for any DG module $M_\idot$
over $A^\hdot$ and every complex $V_\idot$ of $\Q$-vector spaces. On
the level of the stacks $\M(-)$, this tensor product turns
$K^{st}_\idot(A^\hdot)$ into a module over $K^{st}_\idot(\Q) =
\Z[\beta]$. We can now state the main conjecture.

\begin{conj}\label{conj.1}
Assume that $k$ is a ring conatining $\Q$, and assume that a DG
algebra $A^\hdot \Dgalg(k)$ is smooth and proper. Then there exists
a map
$$
c:K^{st}_\idot(A^\hdot) \to HP_\idot(A^\hdot)
$$
such that $c(\beta(\alpha)) = u(c(\alpha))$ for any $\alpha \in
K^{st}_\idot(A^\hdot)$, where $u$ is the periodicity map. The map
$c$ is functorial in $A^\hdot$. Moreover, the induced map
\begin{equation}\label{ov.hp}
K^{st}_\idot(A^\hdot) \otimes_{\Z[\beta]} k[\beta,\beta^{-1}] \to
HP_\idot(A^\hdot)
\end{equation}
is an isomorphism.
\end{conj}

The reason this conjecture is relevant to the present paper is that
the tensor product $K^{st}_\idot(A^\hdot) \otimes k$ by its very
definition has all the structures possessed by the de Rham
cohomology of an algebraic variety. In particular, if $k = \C$,
$K^{st}_\idot(A^\hdot)$ has a canonical real structure.

\begin{conj}\label{conj.2}
  Assume that $K = \C$, and assume given a smooth and proper DG
  algebra $A^\hdot/K$ for which Conjecture~\ref{conj.1} holds.
  Equip $CP_\idot(A^\hdot)$ with the real structure induced from the
  canonical real structure on $K^{st}_\idot(A^\hdot) \otimes K$ by
  the isomorphism ~\ref{ov.hp}. Then for any integer $i$, the
  periodic cyclic homology group $HP_\idot(A^\hdot)$ this real
  structure and the standard Hodge filtration $F^\hdot$ is a pure
  $\R$-Hodge structure of weight $i$.
\end{conj}

The two conjectures above are a slight refinement and/or
reformulation of a conjecture made by B. To\"en \cite{To.lec} with a
reference to A. Bondal and A. Neeman, and described by L. Katzarkov,
M. Kontsevich and T. Pantev in \cite[2.2.6]{KS}.

Apart from the basic case $A^\hdot = k$ of Corollary~\ref{to.c},
the only real evidence for Conjecture~\ref{conj.1} comes from recent
work of Fiedlander and Walker \cite{FW}, where it has
been essentially proved for a DG algebra $A^\hdot$ equivalent to a
smooth projective algebraic variety $X/k$. The definition of
semi-topological $K$-theory used in \cite{FW} is different from
Definition~\ref{K.st}, but it is very close to the homotopy groups
of the group completion of the $E_\infty$-space
$\Top(\overline{\M}(X))$; Lemma~\ref{to.2} should then show that the
two things are the same. Friedlander and Walker also show that their
constructions are compatible with the complex conjugation, so that
Conjecture~\ref{conj.2} then follows by the usual Hodge theory
applied to $X$.

In the general case, as far as I know, both Conjecture~\ref{conj.1}
and Conjecture~\ref{conj.2} are completely open. They are now a
subject of investigation by B. To\"en and A. Blanc.

\bigskip

\subsection*{Acknowledgements.} I have benefited a lot from
discussing this material with A. Beilinson, R. Bezrukavnikov,
A. Bondal, V. Drinfeld, L. Hesselholt, V. Ginzburg, L. Katzarkov,
D. Kazhdan, B. Keller, M. Kontsevich, A. Kuznetsov, N. Markarian,
J.P. May, G. Merzon, D. Orlov, T. Pantev, S.-R. Park, B. To\"en,
M. Verbitsky, G. Vezzosi, and V. Vologodsky; discussions with
M. Kontsevich, on one hand, and B. To\"en and G. Vezzosi, on the
other hand, were particularly invaluable.

\bigskip

\noindent
{\sc
Independent University of Moscow \&
Steklov Math Institute\\
Moscow, USSR}

\bigskip

\noindent
{\em E-mail address\/}: {\tt kaledin@mi.ras.ru}

\end{document}